\documentclass{article}
\usepackage[utf8]{inputenc}

\usepackage{amsmath}
\usepackage{amsfonts}
\usepackage{amsthm}
\usepackage{amssymb}

\usepackage{changepage}

\usepackage[graphicx]{realboxes}

\usepackage{algorithm}
\usepackage{algpseudocode}

\usepackage[pass]{geometry}
\usepackage{caption}

\usepackage[plainpages=false,pdfpagelabels,colorlinks=true,citecolor=blue,hypertexnames=false]{hyperref}

\usepackage{mathtools}

\newtheorem{conj}{Conjecture}

\newtheorem{thm}{Theorem}
\newtheorem{prop}{Proposition}
\newtheorem{lem}{Lemma}

\usepackage{tikz}
\usepackage{tikz-3dplot}

\theoremstyle{definition}
\newtheorem{definition}{Definition}

\DeclareMathOperator{\R}{\mathbb{R}}
\DeclareMathOperator{\N}{\mathbb{N}}
\DeclareMathOperator{\PP}{\mathcal{P}}
\DeclareMathOperator{\QQ}{\mathcal{Q}}
\def\area{\operatorname{Area}}
\def\peri{\operatorname{Peri}}
\def\dia{\operatorname{Dia}}
\def\Rup{\operatorname{Rup}}

\title{An algorithmic approach to Rupert's problem}
\author{Jakob Steininger, Sergey Yurkevich}

\begin{document}

\maketitle

\begin{abstract}
    A polyhedron $\textbf{P} \subset \R^3$ has Rupert's property if a hole can be cut into it, such that a copy of $\textbf{P}$ can pass through this hole. There are several works investigating this property for some specific polyhedra: for example, it is known that all 5 Platonic and 9 out of the 13 Archimedean solids admit Rupert's property. A commonly believed conjecture states that every convex polyhedron is Rupert. We prove that Rupert's problem is algorithmically decidable for polyhedra with algebraic coordinates. We also design a probabilistic algorithm which can efficiently prove that a given polyhedron is Rupert. Using this algorithm we not only confirm this property for the known Platonic and Archimedean solids, but also prove it for one of the remaining Archimedean polyhedra and many others. Moreover, we significantly improve on almost all known Nieuwland numbers and finally conjecture, based on statistical evidence, that the Rhombicosidodecahedron is in fact \textit{not} Rupert.
\end{abstract}

\vspace{-0.4cm}

\section{Introduction}
\begin{figure}
    \centering
        \begin{tikzpicture}[scale=1.5]
    \draw[thick] (0,1.6329932) -- (1.4142,0.8164932) -- (1.4142,-0.8165) -- (0,-1.6329932) -- (-1.4142,-0.8164932) -- (-1.4142,0.8165) -- (0,1.6329932) -- (1.4142,0.8164932);
    \draw[thick] (-1.4142,0.8165) -- (0,0) -- (1.4142,0.8164932);
    \draw[thick] (0,0) -- (0,-1.6329932);
    \draw[dashed, gray] (0,0) -- (0,1.6329932);
    \draw[dashed, gray] (-1.4142,-0.8165) -- (0,0) -- (1.4142,-0.8164932);
    \draw[thick, blue] (-1,1) -- (1,1) -- (1,-1) -- (-1,-1) -- (-1,1) -- (1,1);
    \fill[blue!60, opacity=0.2] (-1,1) -- (1,1) -- (1,-1) -- (-1,-1) -- (-1,1);
    \end{tikzpicture}
    \captionof{figure}{The unit square fits inside the regular hexagon of side length $\sqrt{2/3}$.}
    \label{fig:rupert_cube2}
\end{figure}
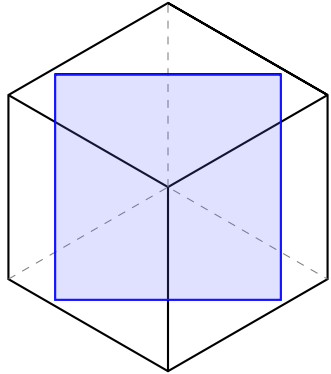

\noindent
Undoubtedly the following fact is surprising when being first encountered with:
\begin{quote}
    \textit{It is possible to cut a hole in the unit cube such that another unit cube can pass through it.}
\end{quote}
Indeed, Prince Rupert of the Rhine won a wager in the 17th century by betting on the validity of this claim. An elegant and simple way to see why this assertion is true is presented in Figure~\ref{fig:rupert_cube2}; indeed, it is easy to verify that the projection of the unit cube in the direction of a main diagonal yields a regular hexagon of side length $\sqrt{2/3}$ and the unit square (a different projection of the cube) fits inside that hexagon. These two observations are already enough to win Rupert's bet, however at the same time they also open a whole world of interesting questions, conjectures and studies.

For instance, a subsequent natural question was investigated by Pieter Nieuwland a century after Prince Rupert's death: 
\begin{quote}
    \textit{How large can the second cube maximal be in order to still fit inside a hole of the unit cube?}
\end{quote}
Quite surprisingly, Nieuwland could show that in terms of this question, the solution presented in Figure~\ref{fig:rupert_cube2} is not optimal. If viewed from a slightly different angle, a ``tunnel'' inside the unit cube can be constructed such that a cube with side length $1.06$ can be moved through it\footnote{The side length of the blue square in Figure~\ref{fig:rupert_cube2} is at most $\sqrt{6} - \sqrt{2} \approx 1.0353$.}. Nieuwland could even find the exact maximal side length of the ``fitting'' cube which turns out to be $3\sqrt{2}/4 \approx 1.06066$ (for a proof of this fact see~\cite{BeGuHuJo21}); later this constant was given the name \textit{Nieuwland's constant}.

Analogously to the cube, Rupert's property can be defined for any polyhedron in $\R^3$. A somewhat imprecise definition of this property is: a polyhedron $\textbf{P} \subset \R^3$ has Rupert's property if a hole (with the shape of a straight tunnel) can be cut into it such that a copy of $\textbf{P}$ can be moved through this hole. In the next section we will first recall a rigorous (but rather non-transparent) definition from~\cite{JeWeYu17} (see Definition~\ref{defi:rupert}) and then an easy end explicit reformulation using projections to $\R^2$ (in the spirit of Figure~\ref{fig:rupert_cube2}). In the same way, the \textit{Nieuwland number} can also be generalized for any polyhedron $\textbf{P}$, see Definition~\ref{def:nieuwlandNumber}.

For a historic overview on these questions we refer to \cite{Schreck50}; for more recent contributions see \cite{Scriba68,JeWeYu17,ChYaZa18,Hoffmann19,Lavau19}. Scriba showed in 1968 that the Tetrahedron and Octahedron have Rupert's property. Half a century later and hence already quite recently, Jerrard, Wetzel and Yuan, the authors of the second paper, built on Scriba's work and investigated Rupert's property of Platonic solids further: they could prove that all five of them are Rupert. Moreover, they also gave lower bounds on Nieuwland numbers for them. One year later Chai, Yuan and Zamfirescu looked at Archimedean solids from ``Rupert's perspective'', showed that 8 out of 13 have Rupert's property and also provided lower bounds for the corresponding Nieuwland numbers. Finally another year later, Hoffmann~\cite{Hoffmann19} and Lavau~\cite{Lavau19} showed in 2019 Rupert's property for the Truncated  tetrahedron, thus enhancing the number to 9 out of 13. Theorem~\ref{thm:truncated_icosidodecahedron} in the present work ``resolves'' the Truncated icosidodecahedron, pushing the number of settled down Archimedean solids to~10.\\

After the submission of this work, a preprint~\cite{ToWi22} by Tonpho and Wichiramala appeared on the internet in which the authors study Rupert's problem in $n$ dimensions and also quote results from a master thesis by Tonpho~\cite{Tonpho18} from~2018. There a relatively similar to parts of the present work (but solely numerical) approach is used to find solutions to Rupert's problem for all Platonic and some Archimedean solids. Even more recently, Fredriksson~\cite{Fredriksson22} built on the ideas from the present work, applied the algorithm for placements of convex polygons from~\cite{AgAmSh98}, and used non-linear optimization techniques like SLSQP and Nelder-Mead in order to obtain new results. Most notably, he was able to prove that the Catalan solids the Triakis tetrahedron and the Pentagonal icositetrahedron have Rupert's property. Finally, an extended abstract of the present work was published in the proceedings of the conference~\href{https://www.issac-conference.org/}{ISSAC~2022}~\cite{StYu22}. There we summarize a selection of our findings and provide explicit conjectural equations for the open Question~\ref{question:3} in~\S\ref{sec:discussion}. \\

\noindent
\textbf{Contribution and structure of the paper}\\ 
In Section~\ref{sec:prelim} we introduce the necessary elementary definitions and concepts. We rigorously define Rupert's property of a polyhedron $\textbf{P}$ and then show that it is equivalent to the existence of a septuple of real numbers satisfying a simple property depending on $\textbf{P}$. In the same section we recall the notion of the Nieuwland number of a polyhedron.

Contrary to the existing methods for proving that a polyhedron has Rupert's property, we present a new algorithmic approach to this problem in Section~\ref{sec:algos}. Roughly speaking, our probabilistic (\textit{Las Vegas} type) algorithm draws pairs of random projections of a given polyhedron and decides whether the chosen directions yield a solution -- if they do not, the algorithm draws another pair, and so on. Moreover, by constructing a deterministic algorithm, we also prove that Rupert's question for most interesting polyhedra is algorithmically decidable. However, we also infer that at least for now this algorithm is only of theoretical value, since it is not yet practical because of its bad complexity. In the same section we explain a simple algorithmic idea which allows to significantly improve on known lower bounds for the Nieuwland numbers. Finally, we also define the concept of the \textit{Rupertness}, measuring the likelihood for finding a solution to Rupert's problem of a (centrally symmetric) polyhedron. 

It turns out that in practice our probabilistic approach finds solutions to Rupert's problem very efficiently: all 5 Platonic and 10 Archimedean solids can be resolved in less than one minute on a regular computer. We present our new explicit results in Section~\ref{sec:results}: we prove that the Truncated icosidodecahedron is Rupert (Theorem~\ref{thm:truncated_icosidodecahedron}), show this property for many Catalan and Johnson solids (Theorems~\ref{thm:catalan_rupert_solids} and \ref{thm:johnson_rupert}), and significantly improve on all known Nieuwland numbers (Table~\ref{tab:Nieuwland_constants}), except the Cube, Octahedron and Cuboctahedron. As mentioned, the Nieuwland constant for the Cube is proven to be optimal and the Nieuwland numbers for the Octahedron and Cuboctahedron  are conjectured to be optimal as well~\cite[p.~91]{JeWeYu17}; our findings support this conjecture. 

In \cite[p.~87]{JeWeYu17} the authors suggest the possible non-existence of ``non-Rupert'' convex polyhedra in $\R^3$ and infer that in any case any such example would be of considerable interest. The authors of \cite{ChYaZa18} go even further and state the following conjecture. 
\begin{conj}[Chai, Yuan, Zamfirescu \cite{JeWeYu17, ChYaZa18}] \label{conj:every_poly_rupert}
Every convex polyhedron has Rupert's property.
\end{conj}
\noindent
Also in Section~\ref{sec:results} we provide statistical evidence for a counter-example to this conjecture (Conjecture~\ref{conj:RCD_non_Rupert}).

Appendix~\ref{sec:appendix} contains our solutions for Platonic, Archimedean and Catalan solids and corresponding lower bounds for the Nieuwland numbers. All these solutions are given in a uniform way in one table. Together with the exact coordinates for the Platonic and Archimedean polyhedra we used (also in the appendix) these solutions can be easily verified by the reader; for the coordinates of Catalan solids we refer to \href{http://dmccooey.com/polyhedra/}{www.dmccooey.com/polyhedra/}. For the reader's convenience, we also provide our source code written in the programming language R and the software Maple: \href{https://github.com/Vog0/RupertProblem/}{www.github.com/Vog0/RupertProblem}.

\section{Preliminaries} \label{sec:prelim}
In order to avoid confusion, let us first collect some elementary definitions. 
\begin{definition}\label{def:def}
The following classical notions we will use throughout the text: 
\begin{itemize}
    \item A \emph{polyhedron}, in this text usually denoted by $\textbf{P}$ or $\textbf{Q}$, is a finite non-degenerate set of points in $\R^3$ in convex position. We denote by $\overline{\textbf{P}}$ the smallest convex set containing all points of $\textbf{P}$ (i.e. including the interior) and by $\textbf{P}^\circ$ its interior.

    \item A \emph{polygon}, usually denoted by $\PP$ or $\QQ$, is a finite set of points in $\R^2$ that not all lie on the same line. Similar to polyhedra, we denote by $\overline{\PP}$ the convex hull of $\PP$ and by $\PP^\circ$ the interior of $\overline{\PP}$.

    \item We call $\Sigma$ the set of isometries of $\R^2$ that do not include reflections, i.e. length preserving mappings from $\R^2$ onto itself not including reflections. It is well-known that any element $\sigma \in \Sigma$ can be represented by a rotation about the origin followed by a translation, and also the other way around. We will let $\sigma\in \Sigma$ act on a set of points in the plane elementwise. Furthermore, we parametrize all translations of $\R^2$ by $T_{x,y}$: $\mathbb{R}^2 \to \mathbb{R}^2$, $$T_{x,y}((a,b)^t)\coloneqq(a+x,b+y)^t.$$
    Similarly, the rotation mapping $R_\alpha$: $\mathbb{R}^2 \to \mathbb{R}^2$ is defined by
    $$R_\alpha((a,b)^t)\coloneqq
 \begin{pmatrix}
    \cos(\alpha) & -\sin(\alpha)  \\ 
    \sin(\alpha) & \cos(\alpha)  \\
\end{pmatrix}
\begin{pmatrix}
    a  \\ 
    b  \\
    \end{pmatrix}
    =
\begin{pmatrix}
    a\cos(\alpha)-b\sin(\alpha)  \\
    a\sin(\alpha)+b\cos(\alpha)  \\
\end{pmatrix}.
    $$    
    Clearly, $T_{x,y}$ translates points in $\R^2$ by the vector $(x,y)^t$ and $R_\alpha$ rotates a point counter-clockwise by an angle $\alpha$ about the origin. 

    \item We say that a polygon $\PP$ \emph{lies inside} a polygon $\QQ$ if $\PP \subset \QQ^\circ$. Moreover, we say that a polygon $\PP$ \emph{fits in} a polygon $\QQ$ if there exists an isometry $\sigma \in \Sigma$ such that $\sigma(\PP)$ lies inside $\QQ$. 
    
    \item A polyhedron $\textbf{P}$ is called \emph{centrally symmetric} with respect to $O_{\textbf{P}} \in \R^3$ if for each $A \in \textbf{P}$, the point $2O_{\textbf{P}}-A$ belongs to $\textbf{P}$. Analogously, a polygon $\PP$ is centrally symmetric about $O_{\PP} \in \R^2$ if $2O_{\PP}-A \in \PP$ for each point $A \in \PP$. A polyhedron or a polygon is called \emph{point symmetric} if it is centrally symmetric with respect to some point.
\end{itemize}
\end{definition}

Usually Rupert's property is explained as follows: a polyhedron $\textbf{P}$ is Rupert if a hole with the shape of a straight tunnel can be cut into it such that a copy of $\textbf{P}$ can be moved through this hole. While this definition explains well why this notion is geometrically intriguing, it is admittedly not quite mathematically precise. A rigorous definition is given for example in \cite{JeWeYu17} and we will state here a slightly reformulated version. First, let us set the notion of a \emph{set with a hole}: we mean a set of points in $\R^3$ whose interior is connected but not simply connected. Given a polyhedron $\textbf{Q}$, we may move it along a straight line in the direction of a vector $v \in \R^3$; taking the convex hull of the union of all these translations we obtain the set $\{ \overline{\textbf{Q}} + tv \in \R^3 \colon t\in \R \}.$ Rupert's property of a polyhedron may be defined as follows (see \cite{JeWeYu17}).

\begin{definition}[Rupert's property]\label{defi:rupert}
    A polyhedron\footnote{Note that, as defined in Definition~\ref{def:def}, a polyhedron in this text is always convex.} $\textbf{P}$ has \emph{Rupert's property} (or $\textbf{P}$ is \emph{Rupert}) if there exists a polyhedron $\textbf{Q}$ of the same shape and size as $\textbf{P}$ and a vector $v \in \R^3$ such that $\overline{\textbf{P}} \setminus \{\overline{\textbf{Q}} + tv \in \R^3 \colon t\in \R \}$
    is a set with  a hole. \emph{Rupert's problem} is the task to decide whether a given polyhedron is Rupert.
\end{definition}

Luckily, the definition of Rupert's property can be reformulated in a much easier criterion on the level of projections to the plane $\R^2$. The idea is that looking from the direction of the vector $v$ in the definition above, we must see the two shadows (normal projections) of the polyhedra $\textbf{P}, {\textbf{Q}}$ as two polygons $\PP,\QQ$, one lying inside the other: $\PP \subset \QQ^\circ$. This is the core of Theorem~1 in \cite{JeWeYu17} and the reason why Figure~\ref{fig:rupert_cube2} in the introduction is a proof that the Cube is Rupert. Now we will make this idea even more explicit.

As we are dealing with projections, we first parametrize the set of all those. We define the mapping $X: [0,2\pi)\times [0,\pi] \to \{x\in \mathbb{R}^3:\lVert x \rVert=1\}$ by
\begin{equation}\label{eq:sphere_parametrization}
X(\theta,\varphi)\coloneqq(\cos \theta \sin \varphi,\sin \theta \sin \varphi, \cos \varphi)^t.
\end{equation}
This gives a way to parametrize the points on the 3-dimensional sphere in terms of two unknowns. It is well-known that drawing $\theta$ uniformly on $(0,2\pi)$, that is $ \theta \sim U(0,2\pi)$, and $\varphi \sim  \arccos(U(-1,1))$ results in a uniformly distributed $X(\theta,\varphi)$ on the unit sphere.

\begin{center}
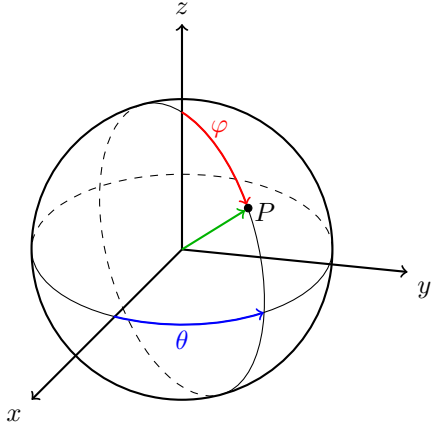

\tikzset{
    partial ellipse/.style args={#1:#2:#3}{
        insert path={+ (#1:#3) arc (#1:#2:#3)}
    }
}

\begin{tikzpicture}
    \draw[thick,->] (0,0) -- (-2,-2) node[anchor=north east]{$x$};
    \draw[thick,->] (0,0) -- (3,-0.3) node[anchor=north west]{$y$};
    \draw[thick,->] (0,0) -- (0,3) node[anchor=south]{$z$};
    \draw[thick] (0,0) circle (2cm);
    \draw (-2,0) arc(180:360:2cm and 1cm);
    \draw[dashed] (2,0) arc(0:180:2cm and 1cm);
    \draw[rotate=105] (-2,0) arc(180:360:2cm and 1cm);
    \draw[rotate=105, dashed] (2,0) arc(0:180:2cm and 1cm);
    
    \draw[thick, blue,->] (0,0) [partial ellipse=243:303:2cm and 1cm];
    \draw[thick,rotate=105, red,<-] (0,0) [partial ellipse=280:332:2cm and 1cm];
    \draw[fill=black] (0.88,0.55) circle (0.05);
    \draw[thick, ->, black!30!green] (0,0)--(0.84,0.52);
    
    \node at (1.1,0.5) {$P$};
    \node at (0.5,1.6) {\color{red} $\varphi$};
    \node at (0,-1.2) {\color{blue} $\theta$};
\end{tikzpicture}
    \captionof{figure}{Meaning of $\theta$ and $\varphi$ in (\ref{eq:sphere_parametrization}) in spherical coordinates.}
    \label{fig:spherical_coords}
\end{center}

\noindent
It follows that a projection onto a plane orthogonal to $X(\theta,\varphi)$ can be given by 
\begin{equation} \label{eq:Mthetaphi}
M_{\theta,\varphi} \coloneqq
\begin{pmatrix} 
-\sin(\theta) & \cos(\theta) & 0 \\ 
-\cos(\theta)\cos(\varphi) & -\sin(\theta)\cos(\varphi) & \sin(\varphi) \\ 
\end{pmatrix}.
\end{equation}

Like the mappings $R_{\alpha},T_{x,y}$ we extend the map $M_{\theta,\varphi}$ to act on sets of points in $\R^3$ elementwise. 
Thus, all parallel projections of the vertices of a  polyhedron $\textbf{P}$ onto $\R^2$ can now be expressed as
$$(T_{x,y}\circ R_\alpha \circ M_{\theta,\varphi}) (\textbf{P}).$$
It follows that an equivalent characterization of Rupert's property for a polyhedron $\textbf{P}$ is the existence of two quintuples of parameters $(x_i,y_i,\alpha_i,\theta_i,\varphi_i)$, $i=1,2$, such that
\[
(T_{x_1,y_1}\circ R_{\alpha_1} \circ M_{\theta_1,\varphi_1}) (\textbf{P}) \subset (T_{x_2,y_2}\circ R_{\alpha_2} \circ M_{\theta_2,\varphi_2}) (\textbf{P})^\circ.
\]
In other words, the polygon on the left-hand side lies inside the polygon on the right-hand side and both polygons are obtained by some orthogonal projection, rotation and translation of $\textbf{P}$. Moreover, this condition can be rewritten as
\[
(R_{-\alpha_2}\circ T_{x_1-x_2,y_1-y_2}\circ R_{\alpha_1} \circ M_{\theta_1,\varphi_1}) (\textbf{P}) \subset M_{\theta_2,\varphi_2} (\textbf{P})^\circ.
\]
Note, that $R_{-\alpha_2}\circ T_{x_1-x_2,y_1-y_2}\circ R_{\alpha_1}$ is an isometry on $\R^2$ and thereby may be expressed as the composition of a single rotation and a translation. Hence, we obtain the following equivalent characterization of Rupert's property.

\begin{prop}\label{prop:Rupert_7parameters}
A polyhedron {\normalfont $\textbf{P}$} satisfies Rupert's property, if and only if there exist 7 parameters $x,y\in \R$, $\alpha,\theta_1,\theta_2 \in [0,2\pi)$ and $\varphi_1,\varphi_2 \in [0,\pi]$ such that
\begin{equation}\label{eq:rupert}
    (T_{x,y}\circ R_{\alpha} \circ M_{\theta_1,\varphi_1}) ({\normalfont\textbf{P}}) \subset  M_{\theta_2,\varphi_2} ({\normalfont\textbf{P}})^\circ.
\end{equation}
\end{prop}

Clearly, any solution of Rupert's property can be translated into these 7 parameters and vice versa. Hence, we will encode a solution to Rupert's problem by a vector $(x,y,\alpha,\theta_1,\theta_2,\varphi_1,\varphi_2) \in \R^7$.

Note that from Proposition~\ref{prop:Rupert_7parameters} it is evident that Rupert's property is a statement about containment of points inside an open set. Since the projection, rotation and translation mappings are continuous, it follows that if there exists a solution to Rupert's problem $v = (x,y,\alpha, \theta_1,\theta_2, \varphi_1,\varphi_2) \in \R^7$, then there must exist an open ball in $\R^7$ around $v$ of solutions. In other words, if a solution exists, then there is a set of solutions with positive (Lebesgue) measure. We will use this observation several times throughout the text.\\

Now let us recall the \emph{Nieuwland number} of a polyhedron. If $\textbf{P}$ is Rupert, then by the consideration above there exists a hole in it in which even a slightly larger copy of $\textbf{P}$ can pass through. Naturally, one may ask for the largest polyhedron similar to $\textbf{P}$ which can also be moved through such a hole. In other words, what is the largest (supremum) number $\nu$ for which there exists a copy of $\textbf{P}$, say $\textbf{Q}$, such that $\nu \textbf{Q}$ can be moved in a straight tunnel through $\textbf{P}$? This number $\nu$ is called the \emph{Nieuwland number} of $\textbf{P}$. It can be defined as in Definition~\ref{defi:rupert}, but in view of the more concrete and useful equivalent formulation in  Proposition~\ref{prop:Rupert_7parameters}, we will define it directly via projections.

\begin{definition}[Nieuwland number]\label{def:nieuwlandNumber}
    The \emph{Nieuwland number} $\nu = \nu(\textbf{P})$ of a polyhedron $\textbf{P}$ is the supremum over all $\mu  \in \R$ for which there exist $x,y\in \R$, $\alpha,\theta_1,\theta_2 \in [0,2\pi)$ and $\varphi_1,\varphi_2 \in [0,\pi]$ such that
    \begin{equation} \label{eq:nieuwland}
    (T_{x,y}\circ R_{\alpha} \circ M_{\theta_1,\varphi_1}) (\mu \textbf{P}) \subset  M_{\theta_2,\varphi_2} (\textbf{P})^\circ.
    \end{equation}
\end{definition}

Clearly, $\textbf{P}$ is Rupert if and only if $\nu(\textbf{P}) > 1$. We note that a typo in~\cite[p.~88]{JeWeYu17} incorrectly states ``$\geq$'' in this inequality. In fact, $\nu(\textbf{P}) \geq 1$ holds for every polyhedron, since if $\mu<1$ in (\ref{eq:nieuwland}), one can take all 7 parameters to be equal to 0 (in other words $\textbf{P} = \textbf{Q}$) and the inclusion holds. As mentioned in the introduction, for a Cube $\textbf{P}$ the Nieuwland number $\nu(\textbf{P})$ is proven to be $3\sqrt{2}/4$. \\

Now we will prove that in the case when the polyhedron is point symmetric, the number of parameters in Proposition~\ref{prop:Rupert_7parameters} can be reduced to 5. This significantly simplifies the algorithms in the next section in the point symmetric case.
\begin{prop}\label{prop:pointsym}
The following two statements hold:
\begin{itemize}
    \item [1)] Let $\PP$ and $\QQ$ be convex polygons which are centrally symmetric around $O_{\PP}$ and $O_{\QQ}$ respectively. Then $\PP$ fits in $\QQ$ if and only there exists a $\sigma\in\Sigma$ such that $\sigma(\PP)$ lies inside $\QQ$ and $\sigma(O_{\PP})=O_{\QQ}$.
    \item [2)] Let ${\normalfont\textbf{P}}$ by a polyhedron that is centrally symmetric about the origin. Then ${\normalfont\textbf{P}}$ satisfies Rupert's property if and only if there are 5 parameters $\alpha\in \R$, $\theta_i,\in [0,2\pi)$ and $\varphi_i \in [0,\pi]$ for $i=1,2$ such that 
    \[
    (R_{\alpha} \circ M_{\theta_1,\varphi_1}) ({\normalfont\textbf{P}}) \subset  M_{\theta_2,\varphi_2} ({\normalfont\textbf{P}})^\circ.
    \]
\end{itemize}
\end{prop}

\begin{proof}
For the first statement it suffices to show that if there exists $\sigma \in \Sigma$ such that $\sigma(\PP)$ lies inside $\QQ$, then there also exists $\sigma' \in\Sigma$ such that $\sigma'(\PP)$ is inside $\QQ$ and $\sigma'(O_{\PP})=O_{\QQ}$. Let $\tau$ be the translation in $\R^2$ which maps $\sigma(O_{\PP})$ to $O_{\QQ}$. We claim that $\sigma' = \tau \circ \sigma$ satisfies the required conditions.

Obviously, $\sigma'(O_{\PP}) = \tau(\sigma(O_{\PP})) = O_{\QQ}$, hence we are left to show that $\sigma'(O_{\PP})$ lies inside $\QQ$. Let $\PP'$ be the reflection of $\sigma(\PP)$ around $O_{\QQ}$. Because $\sigma(\PP)$ is inside $\QQ$ and $\QQ$ is centrally symmetric, $\PP'$ also lies inside $\QQ$. Since $\PP$ is centrally symmetric, it follows that $\PP'$ can be obtained from $\sigma(\PP)$ by a translation. Moreover, $\sigma'(P)$ is given by the aritmetic mean between $\sigma(P)$ and this translation $\PP'$. Now convexity of $\QQ$ implies that $\sigma'(P)$ lies inside $\QQ$.  

We will now prove the second assertion. As both $M_{\theta,\varphi}$ and $R_\alpha$ are linear mappings, one has for any given point $p\in \R^3$ that
\[
(R_\alpha\circ M_{\theta,\varphi})(-p)=-(R_\alpha\circ M_{\theta,\varphi})(p).
\]
Therefore, any pair of antipodal points of the polyhedron is mapped to antipodal points in $\R^2$, resulting in a centrally symmetric polygon about the origin. So the claim follows from the first assertion.
\end{proof}

In order to keep the notation of Proposition~\ref{prop:Rupert_7parameters}, we will encode a solution to Rupert's problem of a point symmetric polyhedron by a 7-dimensional vector $(0,0,\alpha,\theta_1,\theta_2,\varphi_1,\varphi_2) \in \R^7$ as well.

\section{The algorithms}\label{sec:algos}

In this section we present algorithmic ideas for proving or disproving that a given polyhedron $\textbf{P}$ is Rupert. We start by introducing a naive algorithm which searches for a solution to Rupert's problem for a given polyhedron. Then we gradually expand its sophistication and significantly improve the performance. Furthermore, we introduce a method for finding solutions with a high Nieuwland number. We note that all practical algorithms we present are probabilistic of \emph{Las Vegas} type: If a solution is found, it is easy to check (rigorously) its correctness, however the search running time is probabilistic and cannot be known for sure in advance. We explain another viewpoint in \S\ref{sec:deterministic}, where we construct a deterministic algorithm, thus prove that Rupert's problem is algorithmically decidable. However, we also explain that in practice this algorithm is not (yet) useful. Finally, in \S\ref{sec:Rupertness} we introduce the probabilistic concept of the \emph{Rupertness} of a (point symmetric) polyhedron as the likelihood of finding a solution to the corresponding Rupert's problem. 
\subsection{Probabilistic algorithm for solving Rupert's problem} \label{sec:solverupert}
Proposition~\ref{prop:Rupert_7parameters} states that a polyhedron ${\normalfont\textbf{P}}$ satisfies Rupert's property if and only if there are $x,y\in \R$, $\alpha, \theta_i,\in [0,2\pi)$ and $\varphi_i \in [0,\pi]$ for $i=1,2$ such that
\[
(T_{x,y}\circ R_{\alpha} \circ M_{\theta_1,\varphi_1}) ({\normalfont\textbf{P}}) \subset M_{\theta_2,\varphi_2} ({\normalfont\textbf{P}})^\circ.
\]
It seems at first that the two parameters $x$ and $y$ are unbounded. For the first upcoming algorithm it is however necessary to bound all parameters. Hence, we prove the following proposition.
\begin{prop}\label{prop:bound_x,y}
Let ${\normalfont\textbf{P}}$ be a polyhedron containing the origin and let $R \in \R$ be the maximal distance of its vertices to the origin. Assume that a solution to the corresponding Rupert's problem 
\[
(T_{x,y}\circ R_{\alpha} \circ M_{\theta_1,\varphi_1}) ({\normalfont\textbf{P}}) \subset  M_{\theta_2,\varphi_2} ({\normalfont\textbf{P}})^\circ.
\]
is given. Then $|x|,|y|\leq R$.
\end{prop}

\begin{proof}
As $\textbf{P}$ lies inside the ball with radius $R$ centered at the origin, we have
$$M_{\theta_2,\varphi_2}(\textbf{P})\subset \{a\in \R^2:|\!|a|\!|\leq R\}.$$
Since the origin is in the interior of $\textbf{P}$, we have
$$(T_{x,y}\circ R_{\alpha} \circ M_{\theta_1,\varphi_1}) ((0,0,0)^t)=(T_{x,y}\circ R_{\alpha})((0,0)^t)=T_{x,y}((0,0)^t)=(x,y)^t,$$
hence 
$$(x,y)^t\in M_{\theta_2,\varphi_2} (\textbf{P})^\circ \subset \{a\in \R^2:|\!|a|\!|\leq R\}.$$
Therefore $x^2+y^2\leq R^2$ and in particular $|x|,|y|\leq R$.
\end{proof}
Now the interval for each of the 7 parameters in (\ref{eq:rupert}) is bounded and we can create a first version of our probabilistic deciding algorithm.\\

\noindent
\underline{\textsf{Algorithm~1}}\\
\underline{Input}: A polyhedron $\textbf{P}$.\\
\underline{Output}: The solution encoded by $(x,y,\alpha,\theta_1,\theta_2,\varphi_1,\varphi_2) \in \R^7$ if $\textbf{P}$ is Rupert.
\begin{enumerate}
    \item[(1)] Find $R$ like in Proposition~\ref{prop:bound_x,y}. Draw $x$ and $y$ uniformly in $[-R,R]$, $\theta_1$, $\theta_2$ and $\alpha$ uniformly in $[0,2\pi)$, and $\varphi_1$, $\varphi_2$ uniformly in $[0,\pi]$.
    \item[(2)] Construct the two $3 \times 2$ matrices $A$ and $B$ corresponding to the linear maps $R_{\alpha} \circ M_{\theta_1,\varphi_1}$ and $M_{\theta_2,\varphi_2}$. Compute the two projections of $\textbf{P}$ given by ${\mathcal{P}'} \coloneqq T_{x,y}(A \cdot \textbf{P}) =A \cdot \textbf{P} + (x,y)$ and ${\mathcal{Q}'} \coloneqq  B \cdot \textbf{P}$.
    \item[(3)] Find vertices on the convex hulls of ${\mathcal{P}'}$ and ${\mathcal{Q}'}$; denote them by $\mathcal{P}$ and $\mathcal{Q}$.
    \item[(4)] Decide whether $\mathcal{P}$ lies inside of $\mathcal{Q}$ by checking each vertex of $\mathcal{P}$. 
    \item[(5)] If Step~(4) yields a \textsf{True}, return the solution $(x,y,\alpha,\theta_1,\theta_2,\varphi_1,\varphi_2)$. Otherwise, repeat Steps (1)-(5).
\end{enumerate}

\noindent
Here is a pseudocode for this algorithm:

\begin{algorithm}[H]
\begin{algorithmic}[1]
\State $R \gets \sqrt{\max_i(\textbf{P}[i,1]^2+\textbf{P}[i,2]^2+\textbf{P}[i,3]^2)}$
\State $\text{isRupert} \gets \text{False}$
\While{$\text{isRupert} = \text{False}$}
\State Draw $x$ and $y$ uniformly in $[-R,R]$
\State Draw $\theta_1$, $\theta_2$ and $\alpha$ uniformly in $[0,2\pi)$
\State Draw $\varphi_1$ and $\varphi_2$ uniformly in $[0,\pi]$
\State $A \gets R_{\alpha} \circ M_{\theta_1,\varphi_1}$ and $B \gets M_{\theta_2,\varphi_2}$ \Comment $A, B$ are $3 \times 2$ matrices
\State ${\mathcal{P}'}\gets A \cdot \textbf{P} + (x,y)$ and ${\mathcal{Q}'}\gets  B \cdot \textbf{P}$ 
\State $\mathcal{P} \gets \texttt{ConvexHullPoints}({\mathcal{P}'})$ and  $\mathcal{Q} \gets \texttt{ConvexHullPoints}({\mathcal{Q}'})$
\State $n \gets$ \# rows of $\mathcal{P}$ and $m \gets$ \# rows of $\mathcal{Q}$
\For{$i$ from $1$ to $n$}
\State $P \gets \mathcal{P}[i,]$ \Comment{$P\in \R^2$, $i$th row of $\mathcal{P}$ is the $i$th vertex of $\mathcal{P}$}
\If{$P$ is not inside $\mathcal{Q}$}
\State $\text{isRupert} \gets \text{False}$
\State Break the \textbf{For Loop}
\EndIf
\State $\text{isRupert} \gets \text{True}$
\EndFor
\EndWhile
\State \Return $(x,y,\alpha,\theta_1,\theta_2,\varphi_1,\varphi_2)$
\end{algorithmic}
\label{alg:alg0}
\caption{Probabilistic algorithm for deciding whether $\textbf{P}$ is Rupert.\\
\textbf{Input:} Polyhedron $\textbf{P}$ given by an $N \times 3$ matrix for some $N \in \N$.\\
\textbf{Output:} The solution encoded by $(x,y,\alpha,\theta_1,\theta_2,\varphi_1,\varphi_2) \in \R^7$ if $\textbf{P}$ is Rupert.}
\end{algorithm}

Already this very simple algorithm is able to find solutions for many polyhedra. However, it is quite slow, mostly because the 7-dimensional search space for $(x,y,\alpha,\theta_1,\theta_2,\varphi_1,\varphi_2)$ is large. The first and most significant improvement to \textsf{Algorithm~1} is to reduce the parameter search space from $\R^7$ to $\R^4$ by algorithmically finding $x,y$ and $\alpha$ for given $\theta_1,\theta_2,\varphi_1,\varphi_2$. Chazelle~\cite{Chazelle83} found an efficient algorithm for deciding polygon containment under translation and rotation, which we may conveniently apply. Let us call Chazelle's algorithm \textsf{Chazelle}; it takes as input two polygons $\PP$ and $\QQ$ and outputs $(x,y,\alpha)$ such that $(T_{x,y}\circ R_\alpha)(\PP) \subset \QQ$, and \texttt{False} if no such triple exists. 

Exploiting Proposition~\ref{prop:pointsym}, namely that if $\textbf{P}$ is point symmetric then one can choose $x=y=0$, one can significantly simplify the algorithm in the point symmetric case. Namely, one needs to solve the polygon containment problem only under rotation (and not additionally translation) which is a much easier task: we will call this algorithm \textsf{ChazelleR}: its input are two polygons $\PP$ and $\QQ$ and the output is $(0,0,\alpha)$ such that $R_\alpha(\PP) \subset \QQ$, and \texttt{False} if no such $\alpha$ exists. 

We also note that choosing $\theta_1$, $\theta_2$ uniformly in $[0,2\pi)$, and $\varphi_1, \varphi_2$ uniformly in $[0,\pi]$ is slightly unnatural, since this does not give a uniform distribution on the sphere. As explained in \S\ref{sec:prelim}, we will rather draw $\theta_i \sim U(0,2\pi)$ and $\varphi \sim \arccos(U(-1,1))$.
\noindent
We obtain the following improvement to our \textsf{Algorithm~1}:\\

\noindent
\underline{\textsf{Algorithm~2}} (Using \textsf{Chazelle})\\
\underline{Input}: A polyhedron $\textbf{P}$.\\
\underline{Output}: The solution encoded by $(x,y,\alpha,\theta_1,\theta_2,\varphi_1,\varphi_2) \in \R^7$ if $\textbf{P}$ is Rupert.
\begin{enumerate}
    \item[(1)] For each $i \in \{1,2\}$ draw $\theta_i$ uniformly in $[0,2\pi)$, and $\widetilde{\varphi}_i$ uniformly in $[-1,1]$. Set $\varphi_i \coloneqq \arccos(\widetilde{\varphi}_i)$.
    \item[(2)] Construct the two $3 \times 2$ matrices $A$ and $B$ corresponding to the linear maps $M_{\theta_1,\varphi_1}$ and $M_{\theta_2,\varphi_2}$. Compute the two projections of $\textbf{P}$ given by ${\mathcal{P}'} \coloneqq A \cdot \textbf{P}$ and ${\mathcal{Q}'} \coloneqq  B \cdot \textbf{P}$.
    \item[(3)] Find vertices on the convex hulls of ${\mathcal{P}'}$ and ${\mathcal{Q}'}$; denote them by $\mathcal{P}$ and $\mathcal{Q}$.
    \item[(4)] Call \textsf{Chazelle}$(\PP,\QQ)$ (or \textsf{ChazelleR}$(\PP,\QQ)$ if $\textbf{P}$ is point symmetric).
    \item[(5)] If Step~(4) yields a solution $(x,y,\alpha)$, return $(x,y,\alpha,\theta_1,\theta_2,\varphi_1,\varphi_2)$. Otherwise, repeat Steps (1)-(5).
\end{enumerate}

The algorithm above can find solutions to Rupert's problem for many solids in fractions of seconds and is able to solve one of the previously unsolved Archimedean polyhedra (see Theorem~\ref{thm:truncated_icosidodecahedron}). However, we can improve it even further. Analyzing its practical performance, it is clear that the most time consuming part is Step~(4). Heuristically, this is expected because the theoretical complexity of \textsf{Chazelle} is $O(pq^2)$, if $p$ is the number of vertices of $\PP$ and $q$ the number of vertices of $\QQ$~\cite{Chazelle83}, while all other steps in \textsf{Algorithm~2} are at most linear in $N$, the number of vertices of $\textbf{P}$. Therefore, a natural practical improvement to this algorithm would be to discard pairs $(\PP,\QQ)$ already before Step~(4) if it can be algorithmically easily seen that $\PP$ cannot fit inside $\QQ$. Indeed, we can do so by first computing elementary geometric invariants of the polygons. Moreover, these invariants can be computed for a large batch of polygons coming from randomly drawn projections; then we can discard most pairs and need to test only the remaining ones.  

Define \textit{area} and \textit{perimeter} of a polygon in the obvious way and we call the longest line segment inside $\PP$ \textit{the diameter of $\PP$}. Denote the three by $\area(\PP)$, $\peri(\PP)$ and $\dia(\PP)$ respectively. The following easy lemma allows to speed up our search.  

\begin{lem}\label{lem:areadiameterperimeter}
Assume that a convex polygon $\PP$ fits in a polygon $\QQ$ then:
\begin{enumerate}
    \item The area of $\PP$ is smaller than the area of $\QQ$: $\area(\PP) < \area(\QQ)$. \label{enu:1}
    \item The diameter of $\PP$ is smaller than the diameter of $\QQ$: $\dia(\PP) < \dia(\QQ)$.\label{enu:2}
    \item The perimeter of $\PP$ is smaller than the perimeter of $\QQ$: $\peri(\PP) < \peri(\QQ)$.\label{enu:3}
\end{enumerate}
\end{lem}
\begin{proof}
Since area, perimeter and diameter are invariant under translation and rotation, we may assume that $\PP$ not only fits inside $\QQ$ but already lies inside $\QQ$. Then the statements \ref{enu:1} and \ref{enu:2} become evident. Figure~\ref{fig:perimeter} proves part \ref{enu:3} of the lemma. Note that convexity is important only for this part.
\end{proof}
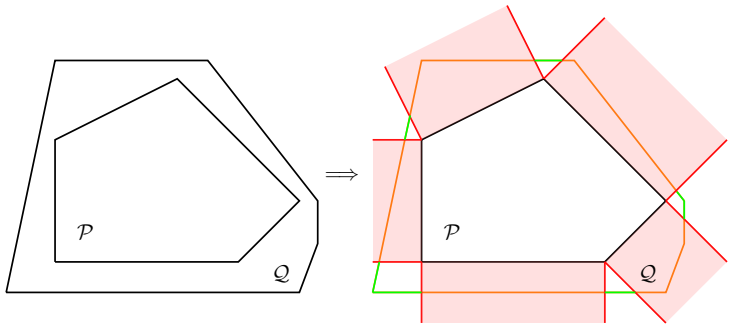
\begin{figure}[b]
    \centering
    \resizebox{0.8\textwidth}{!}{
\begin{tikzpicture}
\draw[thick] (0,0) -- (3,0) -- (4,1) -- (2,3) -- (0,2) -- (0,0);

\draw[thick] (-0.8,-0.5) -- (4,-0.5) -- (4.3,0.3) -- (4.3,1) -- (2.5,3.3) -- (0,3.3) -- (-0.8,-0.5);

\draw[thick] (6,0) -- (9,0) -- (10,1) -- (8,3) -- (6,2) -- (6,0);

\draw[thick, orange] (5.2,-0.5) -- (10,-0.5) -- (10.3,0.3) -- (10.3,1) -- (8.5,3.3) -- (6,3.3) -- (5.2,-0.5);

\draw[thick, red] (6,0) -- (6,-1);
\draw[thick, red] (9,0) -- (9,-1);
\fill[red!60, opacity=0.2] (6,0) rectangle (9,-1);

\draw[thick, red] (9,0) -- (10,-1);
\draw[thick, red] (10,1) -- (11,0);
\fill[red!60, opacity=0.2] (9,0)--(10,-1)--(11,0)--(10,1);

\draw[thick, red] (10,1) -- (11,2);
\draw[thick, red] (8,3) -- (9,4);
\fill[red!60, opacity=0.2] (10,1)--(11,2)--(9,4)--(8,3);

\draw[thick, red] (8,3) -- (7.4,4.2);
\draw[thick, red] (6,2) -- (5.4,3.2);
\fill[red!60, opacity=0.2] (8,3)--(7.4,4.2)--(5.4,3.2)--(6,2);

\draw[thick, red] (6,2) -- (5.2,2);
\draw[thick, red] (6,0) -- (5.2,0);
\fill[red!60, opacity=0.2] (6,2)--(5.2,2)--(5.2,0)--(6,0);

\draw[thick, green] (5.2,-0.5) -- (6,-0.5);
\draw[thick, green] (9,-0.5) -- (9.5,-0.5);
\draw[thick, green] (10.3,0.7) -- (10.3,1);
\draw[thick, green] (10.3,1) -- (10.18,1.15);
\draw[thick, green] (8.3,3.3) -- (7.85,3.3);
\draw[thick, green] (5.815,2.4) -- (5.725,2);
\draw[thick, green] (5.31,0) -- (5.2,-0.5);

\node at (0.5,0.5) {$\PP$};
\node at (3.7,-0.2) {$\QQ$};
\node at (6.5,0.5) {$\PP$};
\node at (9.7,-0.2) {$\QQ$};
\node at (4.7,1.4) {$\Longrightarrow$};
\end{tikzpicture}}
    \caption{Proof that if convex $\PP$ lies inside $\QQ$ then $\peri(\PP) < \peri(\QQ)$.}
    \label{fig:perimeter}
\end{figure}    

Obviously the perimeter of a polygon can be computed in linear time depending on the number of vertices. The Shoelace formula allows for the same complexity for the area. The method of \textit{rotating calipers} allows to compute the diameter of a (convex) polygon in linear time as well~\cite[Cor.~3.1]{Shamos78}. We obtain the following efficient algorithm.\\

\noindent
\underline{\textsf{Algorithm~3}} (Using \textsf{Chazelle} and Lemma~\ref{lem:areadiameterperimeter})\\
\underline{Input}: A polyhedron $\textbf{P}$, a batch size $M \in \N$.\\
\underline{Output}: The solution encoded by $(x,y,\alpha,\theta_1,\theta_2,\varphi_1,\varphi_2) \in \R^7$ if $\textbf{P}$ is Rupert.
\begin{enumerate}
    \item[(1)] For each $j \in \{1,\dots,M\}$ draw $\theta_j$ uniformly in $[0,2\pi)$, and $\widetilde{\varphi}_j$ uniformly in $[-1,1]$. Set $\varphi_j \coloneqq \arccos(\widetilde{\varphi}_j)$.
    \item[(2)] For each $j \in \{1,\dots,M\}$ construct the $3 \times 2$ matrix $A_j$ corresponding to the linear map $M_{\theta_j,\varphi_j}$. Compute the projection of $\textbf{P}$ given by ${\mathcal{P}'_j} \coloneqq A_j \cdot \textbf{P}$. Find the vertices on the convex hull of ${\mathcal{P}'_j}$ and denote them by $\mathcal{P}_j$. Compute and store: $\area(\PP_j)$, $\peri(\PP_j)$ and  $\dia(\PP_j)$.
    \item[(3)] For each $j \in \{1,\dots,M\}$ and $k \in \{ 1,\dots,M\}$ such that $k \neq j$: if $\area(\PP_j) < \area(\PP_k)$ and $\peri(\PP_j) < \peri(\PP_k)$ and  $\dia(\PP_j) < \dia(\PP_k)$ then call \textsf{Chazelle}$(\PP_j,\PP_k)$ (or \textsf{ChazelleR}$(\PP_j,\PP_k)$ if $\textbf{P}$ is point symmetric).
    \item[(4)] If for some pair $(j,k)$ Step~(3) yields a solution $(x,y,\alpha)$, then return $(x,y,\alpha,\theta_j,\theta_j,\varphi_k,\varphi_k)$. Otherwise, repeat steps (1)-(4).
\end{enumerate}
We ran our implementations on all 5~Platonic, 13 Archimedean, 13~Catalan and 92~Johnson polyhedra. The results are presented in Section~\ref{sec:results}.\\

\subsection{Finding and improving Nieuwland's numbers}\label{sec:findimproveNieuwland}

In this short section, we briefly explain an algorithmic method which yields lower bounds on Nieuwland numbers of polyhedra and a simple procedure for finding ``good'' solutions to Rupert's problem. Recall from Definition \ref{def:nieuwlandNumber} that if 
there exist $x,y\in \R$, $\alpha,\theta_1,\theta_2 \in [0,2\pi)$, $\varphi_1,\varphi_2 \in [0,\pi]$ and $\mu \geq 1$ such that
    \begin{equation} \label{eq:nieuwlandV2}
    (T_{x,y}\circ R_{\alpha} \circ M_{\theta_1,\varphi_1}) (\mu \textbf{P}) \subset  M_{\theta_2,\varphi_2} (\textbf{P})^\circ
    \end{equation}
then the Nieuwland number of $\textbf{P}$ is greater than $\mu$, that is $\nu(\textbf{P}) > \mu$. We will say that the Nieuwland number of \emph{a solution} $(x,y,\alpha,\theta_1,\varphi_1,\theta_2,\varphi_2)$ to Rupert's problem of some polyhedron is the largest real number $\mu$ such that (\ref{eq:nieuwlandV2}) holds. Naturally we will say that a solution is better than another if it has a larger Nieuwland number. Clearly, the Nieuwland number of a solution to Rupert's problem for some polyhedron $\textbf{P}$ gives a lower bound on $\nu(\textbf{P})$.

In the previous section we introduced algorithms for finding solutions to Rupert's problem, i.e. finding solutions to $\eqref{eq:nieuwlandV2}$ with $\mu>1$. Given such a solution $(x,y,\alpha,\theta_1,\varphi_1,\theta_2,\varphi_2)$, it is easy to efficiently  find (numerically) an approximation with any given precision for its Nieuwland number using a binary search method: Given \textbf{P}, $v = (x,y,\alpha,\theta_1,\varphi_1,\theta_2,\varphi_2)$ and some $\mu$ it is easy to check whether (\ref{eq:nieuwlandV2}) holds, therefore one can search for the correct $\mu$ by constantly halving the interval which it contains. We will denote this procedure $\mu(v,\textbf{P})$.

Since we are also interested in ``optimal'' solutions to Rupert's problem, i.e. solutions with maximal Nieuwland number, we will briefly explain a procedure to improve a found solution. The idea is simple: starting with a solution $v = (x,y,\alpha,\theta_1,\theta_2,\varphi_1,\varphi_2) \in \R^7$, we first compute $\mu(v, \textbf{P})$ and then perturb all parameters by small random numbers $r_1,\dots,r_7$: $v_i'=v_i+r_i$. If by chance we find a better solution, i.e. if $\mu(v',\textbf{P})>\mu(v,\textbf{P})$, we continue with $v'$, otherwise we choose another random vector $(r_1,\dots,r_7)$. In practice, the numbers $r_i$ are drawn uniformly from some small intervals around 0, which are steadily narrowed down if no improvement was observed for a long time. Also, in order to avoid convergence to ``local minima'', one should run this procedure on several different initial solutions. Similarly to the considerations before, if $\textbf{P}$ is point symmetric, we can choose $x=y=0$ and then draw only $(r_1,\dots,r_5) \in \R^5$.

We note that this method is indeed rather naive and probably may be improved easily. For example, in the very recent work~\cite{Fredriksson22} Fredriksson uses non-linear optimization methods like SLSQP and Nelder-Mead to find numerically optimal solutions for Rupert's problem. Still, in practice we observed that our approach performs quite well. For example, after less than one minute of computational time on a regular computer, we found a solution to Rupert's problem for the Cube and improved it to have Nieuwland's number of 1.06058. As mentioned in the introduction, it is known that the optimal solution for the Cube has Nieuwland's number $3\sqrt{2}/4 \approx 1.06066$.

We will present our results on improved lower bounds for Nieuwland numbers for various solids in \S\ref{sec:lowerNieuwland}.  \\

\subsection{Deterministic algorithm}\label{sec:deterministic}
In this section we will design a deterministic algorithm for deciding whether a given polyhedron satisfies Rupert's property. The main idea is to transform the problem into systems of polynomial inequalities and consequently into the decidability problem of emptiness of semi-algebraic sets.

The first step towards this algorithm is to develop an algebraic formulation for expressing the containment of a point $B\in \R^2$ in the convex hull formed by some points $A_1,\dots,A_n\in \R^2$.

\begin{lem}\label{prop:det>0}
Let $A_1,\dots,A_n\in \R^2$ be the vertices of a convex polygon ordered in counter-clockwise direction and $B\in \R^2$ a point strictly inside this polygon. Set $A_{n+1}\coloneqq A_1$. Then $\det(A_{i}-B,A_{i+1}-B)>0$ for $i=1,\dots,n$.
\end{lem}

\begin{proof}
As $B$ lies inside the described polygon, the oriented angles $\sphericalangle A_iBA_{i+1}$ lie in the open interval $(0,\pi)$. This  implies 
\[
\det(A_i-B,A_{i+1}-B)=\underbrace{|\!|A_i-B|\!|}_{>0}\cdot \underbrace{|\!|A_{i+1}-B|\!|}_{>0} \cdot \underbrace{\sin(\sphericalangle A_{i}BA_{i+1})}_{>0}>0. \qedhere
\]
\end{proof}

\begin{lem}\label{prop:Binconvexhull}
Let $A_1,\dots,A_n,B\in \R^2$, set $A_{n+1}\coloneqq A_1$ and assume that for $i=1,\dots,n$ it holds that $\det(A_{i}-B,A_{i+1}-B) >0$. Then $B$ lies strictly inside the convex hull spanned by $A_1,\dots,A_n$. 
\end{lem}

\begin{proof}
    Assume that $B$ is not inside the interior of the convex hull. By the continuity of the determinant, if $B$ lies exactly on the border, there exists a $B' \in \R^2$ outside the convex hull, still satisfying all (strict) inequalities. So we can assume that $B$ lies outside the convex hull. Then there exists a $v\in \R^2$ such that all $A_i$ strictly lie on the same side of the line ${\{B+vt \colon t\in \R\}\subset \R^2}$. Hence, there is a $w\in \R^2$ perpendicular to $v$ such that every $A_i$ can be written as
    $A_i=B+t_iv+s_iw$, with $t_i\in \R$ and $s_i \in \R^+$. Let $U = (v,w) \in \R^{2 \times 2}$; it follows that 
    \[
    (A_{i} - B, A_{i+1} - B) = U \cdot
    \begin{pmatrix}
        t_{i} & t_{i+1}\\
        s_{i} & s_{i+1}
    \end{pmatrix}.
    \]
    Taking the determinant, we find $0<\det(A_{i}-B,A_{i+1}-B) = \det(U)(t_{i}s_{i+1}-s_it_{i+1})$.
    Dividing by $s_is_{i+1}>0$ yields
    $ 0<\det(U)\left({t_i}/{s_i}-{t_{i+1}}/{s_{i+1}}\right). $
    Finally, summing over all these inequalities gives the desired contradiction.
\end{proof}

Let $\textbf{P}=\{P_1,\dots,P_n\}$ be a convex polyhedron with $n$ enumerated vertices and further let a parallel projection $\PP=(T_{x,y}\circ R_{\alpha} \circ M_{\theta,\varphi}) (\textbf{P})$ of the polyhedron be given. Only a subset of the projected $P_i$ lie on the boundary of $\overline{\PP}$. Let those be $P_{s_1},\dots,P_{s_k}$ ordered in counter-clockwise direction as they appear along the boundary. We call the cycle $s=(s_1,\dots,s_k)$ the \emph{silhouette} of the polyhedron under the projection.

Note that two projections $(T_{x,y}\circ R_{\alpha} \circ M_{\theta,\varphi})(\textbf{P})$ and $(M_{\theta,\varphi})(\textbf{P})$ always have the same silhouette, as translations and rotations do not influence which points of a polygon are on its boundary.

Further, we define $S_n$ to be the set of all non-empty cycles of any (non-empty) subset of the numbers from $1$ to $n$. For instance, we have $S_3=\{(1),(2),(3),(1,2),\allowbreak (2,3),(1,3),(1,2,3),(1,3,2)\}$. 

Clearly, the silhouette of a polyhedron with $n$ enumerated vertices under any projection is an element of $S_n$. The following argument bounds $|S_n|$ from above: Denote by $k$ the length of a cycle and recall that there are $(k-1)!$ cycles of $k$ elements. Hence we have
\begin{equation}\label{eq:Sn<en!}
|S_n|=\sum_{k=1}^n \binom{n}{k}(k-1)!=\sum_{k=1}^n\frac{n!}{k(n-k)!}<n!\sum_{k=1}^n\frac{1}{(n-k)!}<e \cdot n!,
\end{equation}
where $e \approx 2.72$ is Euler's number.

\begin{thm}\label{thm:deterministic}
Let ${\normalfont\textbf{P}}$ be a convex polyhedron with $n$ vertices having integer coordinates, whose absolute value is bounded by $m$. There exists a deterministic algorithm with running time $(\log(m)\cdot n)^{O(1)} \cdot  n!$ deciding whether ${\normalfont\textbf{P}}$ is Rupert and finding a solution if it exists.
\end{thm}

\begin{proof}
We start by enumerating the vertices of the polyhedron $\textbf{P}=\{P_1,\dots,P_n\}$. The algorithm we will present can decide whether there exists a solution to Rupert's problem
\[
    (T_{x,y}\circ R_{\alpha} \circ M_{\theta_1,\varphi_1}) (\textbf{P}) \subset M_{\theta_2,\varphi_2} (\textbf{P})^\circ
\]
for any possible silhouette $s\in S_n$ of the projection on the right-hand side. Then the full algorithm will run over all elements of $S_n$. 

Let $x,y,\alpha,\theta_1,\varphi_1,\theta_2,\varphi_2$ be variables. Given a silhouette $s=(s_1,s_2,\dots,s_k)$, let $Q_i \coloneqq M_{\theta_2,\varphi_2} (\textbf{P}_{s_i})$ and $P_j \coloneqq (T_{x,y}\circ R_{\alpha} \circ M_{\theta_1,\varphi_1}) (\textbf{P}_j)$ for $i=1,\dots,k$ and $j=1,\dots,n$. We also set $Q_{k+1} \coloneqq Q_1$. In other words, $Q_1,\dots,Q_k$ denote the vertices on the boundary of $M_{\theta_2,\varphi_2} (\textbf{P})^\circ$ given a solution with silhouette $s$. Recall that by definition the vertices $Q_{s_1},\dots,Q_{s_k}$ are in ordered in counter-clockwise direction. We define the system of $kn$ inequalities in the seven unknowns $x,y,\alpha,\theta_1,\varphi_1,\theta_2,\varphi_2$:
\begin{equation} \label{eq:detZij>0}
\det({Q_{s_{i}}-P_j,Q_{s_{i+1}}-P_j})>0 \qquad \text{$j=1,\dots,n$ and $i=1,\dots,k$}.
\end{equation}
Now there are two important observations:
\begin{enumerate}
    \item If this system has a solution $(x,y,\alpha,\theta_1,\theta_2,\varphi_1,\varphi_2) \in \R^7$, then by Lemma~\ref{prop:Binconvexhull} all $P_j$ lie in the interior of the convex hull of the $Q_{s_i}$. Therefore this septuple gives a solution to Rupert's problem for $\textbf{P}$ (not necessarily for the silhouette~$s$). 
    \item If the system (\ref{eq:detZij>0}) does not have a solution, then there does not exist a solution to Rupert's problem with the silhouette $s$. In other words, if Rupert's problem for $\textbf{P}$ has a solution $(x,y,\alpha,\theta_1,\theta_2,\varphi_1,\varphi_2)$ for some $s$, then (\ref{eq:detZij>0}) must hold at this point. Since in the definition of silhouette, the vertices are required to be ordered in counter-clockwise direction, we can apply Lemma~\ref{prop:det>0} and the observation follows.
\end{enumerate}
Therefore solving the system (\ref{eq:detZij>0}) is of crucial importance. Denote by $Z_{j,i}$ the matrices $({Q_{s_{i}}-P_j,Q_{s_{i+1}}-P_j})$, i.e. write (\ref{eq:detZij>0}) as $\det(Z_{j,i})>0$. 

Now we would like to employ algorithms for deciding existence of solutions to systems of polynomial inequalities, but the system (\ref{eq:detZij>0}) involves trigonometric functions. However, it is also easy to see that (\ref{eq:detZij>0}) is a polynomial system in the ``variables'' $x,y, \sin(\alpha), \cos(\alpha)$, $ \sin(\theta_i), \cos(\theta_i), \sin(\varphi_i), \cos(\varphi_i)$, $i=1,2$. Henceforth, we shall apply the following rational parametrization of the circle:
\begin{align*}
    f:\mathbb{R} & \to \mathbb{R}^2\\
    t & \mapsto \left(\frac{1-t^2}{1+t^2},\frac{2t}{1+t^2}\right).
\end{align*} 
It is well-known that not only $||f(t)||=1$ for all $t$, but also that $f$ is a bijection between $\mathbb{R}$ and $S^2\setminus \{(-1,0)\}$. We will substitute the variables $\alpha,\theta_i,\varphi_i$, $i=1,2$ with the variables $a,b_1,b_2,c_1,c_2 \in \mathbb{R}$ by
\begin{align*}
\big(\cos(\alpha),\sin(\alpha)\big)&\eqqcolon f(a),\\
\big(\cos(\theta_i),\sin(\theta_i)\big)&\eqqcolon f(b_i),\\
\big(\cos(\varphi_i),\sin(\varphi_i)\big)&\eqqcolon f(c_i),
\end{align*}
for $i=1,2$. Now all entries of $Z_{j,i}$ are rational functions and so are also the inequalities $\det(Z_{j,i})>0$. Next, for each $i$ and $j$ we define the matrix $\widetilde{Z_{j,i}}$ as the matrix $Z_{j,i}$ multiplied by  $(1+a^2)(1+b_1^2)(1+b_2^2)(1+c_1^2)(1+c_2^2)>0$. Each entry of the matrix $\widetilde{Z_{j,i}}$ is a polynomial in $x,y,a,b_1,b_2,c_1,c_2$.

Note that the determinants of $Z_{j,i}$ and $\widetilde{Z_{j,i}}$ have the same sign by the linearity of the determinant. Therefore, the system (\ref{eq:detZij>0}) is equivalent to the system $\det(\widetilde{Z_{j,i}})>0$. Expanding $\det(\widetilde{Z_{j,i}})$ shows that its coefficients are bounded by $O(m^2)$ and the polynomials have a total degree of at most {22}.

Therefore, we are left with a system $\det(\widetilde{Z_{j,i}})>0$ consisting of $kn$ polynomial inequalities in 7 variables, each of them having total degree of at most $22$ and integer coefficients bounded in absolute value by $O(m^2)$. According to \cite{GrVo88}, this system can be solved in a complexity that is polynomial in $\log(m^2) (nk\cdot 22)^{7^2}$, i.e. polynomial in $\log(m)(nk)^{7^2}$. Using $k\leq n$ the complexity simplifies to $(\log(m)n)^{O(1)}$.

Finally, in the worst case, we need to solve such a system for every possible cycle in $S_n$ of possible silhouettes, so using the observation (\ref{eq:Sn<en!}), we get the total upper bound for the running time complexity: $(\log(m) \cdot n)^{O(1)} \cdot  n!$. 
\end{proof}
\noindent
The described algorithm can be summarized as follows:\\

\noindent
\underline{\textsf{Algorithm~4}}\\
\underline{Input}: A polyhedron $\textbf{P}=\{P_1,\dots, P_n\} \subseteq \mathbb{Z}^3$.\\
\underline{Output}: The solution encoded by $(x,y,\alpha,\theta_1,\theta_2,\varphi_1,\varphi_2) \in \R^7$ if $\textbf{P}$ is Rupert.\\

\noindent For every possible silhouette $s=(s_1,\dots,s_k)\in S_n$:

\begin{enumerate}
    \item [(1)] Define the system of inequalities $\det({Q_{s_{i}}-P_j,Q_{s_{i+1}}-P_j})>0$ for $j=1,\dots,n$ and $i=1,\dots,k$,
    where $Q_i \coloneqq M_{\theta_2,\varphi_2} (\textbf{P}_{s_i})$ and $P_j \coloneqq (T_{x,y}\circ R_{\alpha} \circ M_{\theta_1,\varphi_1}) (\textbf{P}_j)$ as well as $Q_{k+1} \coloneqq Q_1$.
    \item [(2)] Substitute the variables $\alpha,\theta_i,\varphi_i$ with $a,b_i,c_i$, $i=1,2$, using the above defined function $f$. This yields a system of rational inequalities.
    \item [(3)] Multiply each inequality by {$\big((1+a^2)(1+b_1^2)(1+b_2^2)(1+c_1^2)(1+c_2^2)\big)^2$}, to get a system of polynomial inequalities with integer coefficients.
    \item [(4)] Search for a solution using the algorithm described in \cite{GrVo88}.
    \item [(5)] If (4) yielded a solution: Transform the found solution back to the original variables $(x,y,\alpha,\theta_1,\theta_2,\varphi_1,\varphi_2) \in \R^7$ using $f^{-1}$. Break the loop and return this septuple as a solution to Rupert's Problem.
\end{enumerate}

Note that Theorem~\ref{thm:deterministic} above can easily be extended to incorporate polyhedra having rational coordinates, as these can be stretched by the least common multiple of the denominators of $\textbf{P}$ in order to have integer coefficients. Moreover, if the coordinates of the polyhedron are not rational but algebraic numbers (like for most Platonic and Archimedean solids) the algorithm above can be adapted as well. The trick is to add to the system of inequalities (\ref{eq:detZij>0}) new variables and equations given by minimal polynomials encoding these coordinates.

We remark that the bound $O(n!)$ on the possible number of silhouettes is very pessimistic. For example, up to (isomorphic) permutations, the Cube has essentially only one silhouette, while $8!$ is quite huge. We are confident that by a closer inspection one can show that the number of possible silhouettes actually growths polynomially in $n$ and for regular polyhedra, like the Platonic or Archimedean solids, is quite small. In practice, however, this does not change much, because the complexity to solve already one single system of inequalities corresponding to a silhouette seems to be infeasible (we will address this issue in \S\ref{sec:results} more explicitly). Therefore, any possible way to reduce the number of silhouettes one needs to check still leads to an algorithm that is unlikely to determine the existence or non-existence of a solution for a non-trivial polyhedron. So we conclude that, at least for now, the described deterministic algorithm is only of theoretical value.

\subsection{Rupertness} \label{sec:Rupertness}

In this section we will quantify the likelihood of finding a solution to Rupert's problem by a randomly chosen projection. For a given polyhedron $\textbf{P}$ we will define the Rupertness $\Rup(\textbf{P})$ as the probability that two random projections of it yield a solution to Rupert's problem. We already discussed that point symmetry is advantageous in general for proving Rupert's property, as it decreases the search space from $\R^7$ to $\R^5$. Keeping that in mind, we will only focus on comparing point symmetric polyhedra and define Rupertness only in this setting:
\begin{definition} \label{def:rupertness}
Let $\textbf{P}$ be a centrally symmetric polyhedron. The Rupertness of~$\textbf{P}$, denoted $\Rup(\textbf{P})$, is the probability that two  uniformly chosen projections $M_{\theta_1,\varphi_1}(\textbf{P}), \allowbreak M_{\theta_2,\varphi_2}(\textbf{P})$ can be extended to a solution of Rupert's problem for $\textbf{P}$, i.e. there exists some $\alpha \in [0,\pi)$ such that 
$(R_{\alpha} \circ M_{\theta_1,\varphi_1}) (\textbf{P}) \subset M_{\theta_2,\varphi_2} (\textbf{P})^\circ.$
\end{definition}
\bigskip

Note that, naturally, in this definition we draw $\theta_i$ and $\varphi_i$ ($i=1,2$) not uniformly on the intervals $[0,2\pi)$ and $[0,\pi]$ but in a way such that the projections are uniformly distributed on the sphere. As mentioned in \S\ref{sec:prelim} this can be modeled by choosing $\theta_i \sim U(0,2\pi)$ uniformly and $\varphi \sim \arccos(U(-1,1))$.

As observed in Section~\ref{sec:prelim}, if $\textbf{P}$ is Rupert then there must already exist a set of solutions with positive measure. Therefore, a point symmetric polyhedron $\textbf{P}$ is Rupert if and only if $\Rup(\textbf{P}) >0$. This also proves that if a solution to Rupert's problem of a polyhedron exists, \textsf{Algorithm~3} will find it eventually.

As we will elaborate in \S\ref{sec:results}, our algorithms can solve all Archimedean polyhedra except three: The Rhombicosidodecahedron (RID in short), Snub cube and Snub dodecahedron. Hence, the RID is the only remaining point symmetric Archimedean polyhedron, for which Rupert's problem is open.
Our main application of the notion of Rupertness is to statistically show that the RID is significantly different from the solved Archimedean polyhedra. 

Using the algorithms from \S\ref{sec:solverupert} and elementary statistics, we can estimate confidence intervals of Rupertness for various solids. For example, if 1000 random pairs of projections of the Cube gave 65 solutions, the probability estimate would be 6.5\% and since this can be viewed as a Bernoulli experiment, one can also calculate the $1-\alpha$ confidence interval $6.5\% \pm \epsilon$ for this probability for any $\alpha \in (0,1)$. More precisely, if $n$ random pairs of projections $M_{\theta_1,\varphi_1}(\textbf{P}),M_{\theta_2,\varphi_2}(\textbf{P})$ gave $k>0$ solutions then the Clopper-Pearson formula implies that the $1-\alpha$ confidence interval for the underlying probability is given by $(S_{\min},S_{\max} )$, where
\begin{align}\label{eq:Clopper-Pearson}
    S_{\min} &=  \left(1 + \frac{n-k+1}{k \cdot F(\alpha/2 \, ;\,2k,2(n-k+1))} \right)^{-1},\\ 
    S_{\max} &= \left(1 + \frac{n-k}{(k+1) \cdot F(1-\alpha/2\,;\,2(k+1),2(n-k))} \right)^{-1},\nonumber
\end{align}
and $F(q \,;\, d_1, d_2)$ is the $q$ quantile of the $F$-distribution with $d_1$ and $d_2$ degrees of freedom. In the case $k=0$, the probability is between $0$ and $1-\sqrt[n]{\alpha/2}$ with a certainty of $1-\alpha$. 

\section{Explicit results} \label{sec:results}
In this section we collect the explicit results of our work. We prove Rupert's property for a tenth Archimedean solid, show that most Catalan and Johnson solids are Rupert, improve on almost all known Nieuwland numbers and estimate the Rupertness of all point symmetric Platonic and Archimedean polyhedra. The solutions described below in Theorem~\ref{thm:truncated_icosidodecahedron} and \ref{thm:catalan_rupert_solids} as well as in \S\ref{sec:lowerNieuwland} are found using the probabilistic and numerical algorithms from the previous section in the programming language R and then verified with rigorous bounds in Maple. 

\subsection{Rupert solids}
We start by resolving a new Archimedean solid:

\begin{thm}\label{thm:truncated_icosidodecahedron}
The Truncated icosidodecahedron has Rupert's property.
\end{thm}

\begin{figure}
    \centering
    \includegraphics[width=7cm]{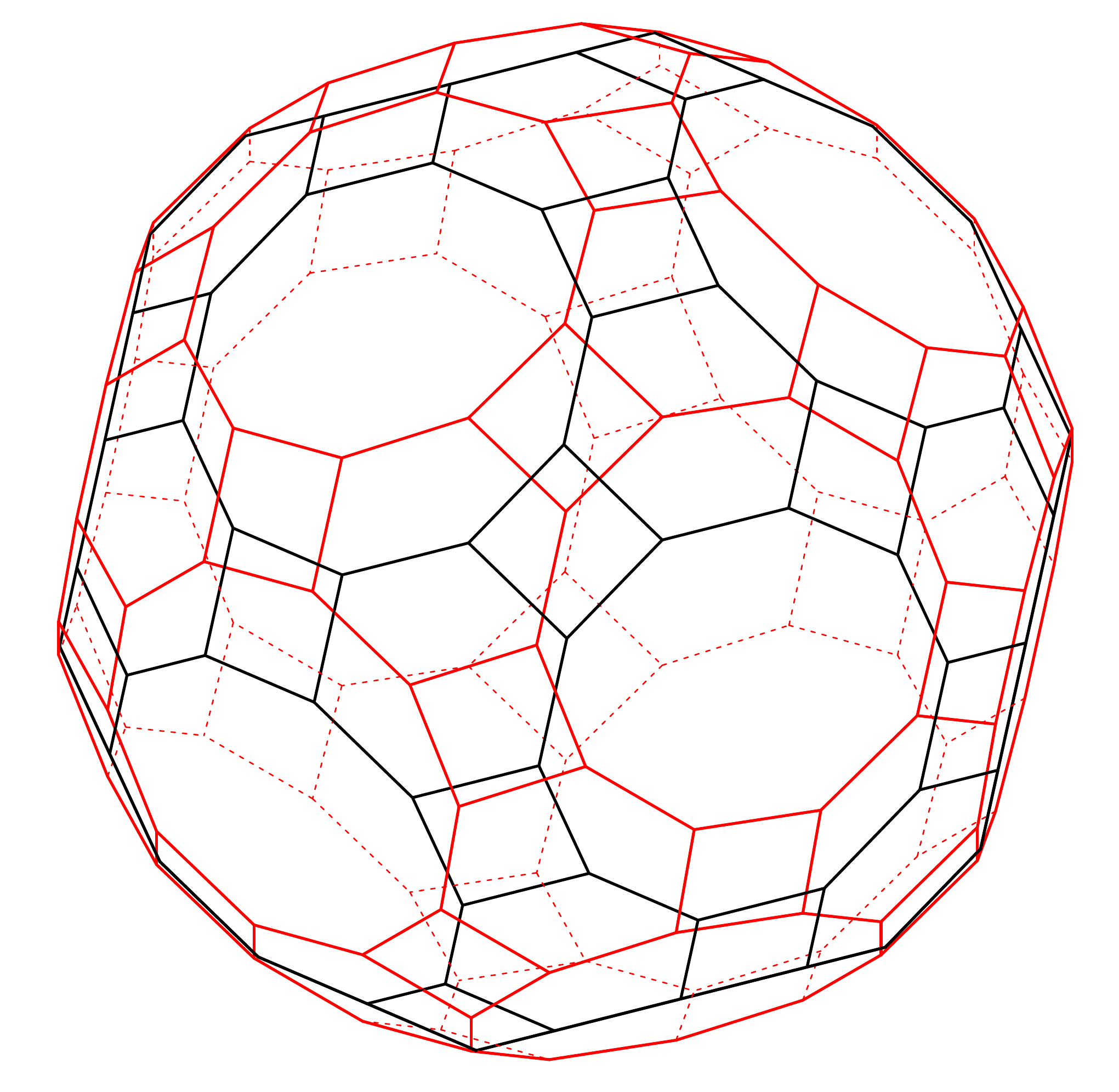}
    \caption{The Truncated icosidodecahedron is Rupert.}
    \label{fig:truncated_icosidodecahedron}
\end{figure}

\begin{proof}
Since this polyhedron is centrally symmetric, we can set $x=y=0$ by Proposition~\ref{prop:pointsym}. So we just need to find the five parameters $ \alpha, \theta_1,\theta_2, \varphi_1, \varphi_2$ as in Proposition~\ref{prop:Rupert_7parameters}. They can be found quickly by applying \textsf{Algorithm~3} to the list of coordinates of the vertices of the Truncated icosidodecahedron (see Table~\ref{tab:Coordinates}). Here is an improved solution (after application of the method described in \S\ref{sec:findimproveNieuwland}): 
\begin{align*}
    & \qquad \qquad \quad \alpha = 0.43584,\\
    & \theta_1 = 2.77685, \quad \quad  \theta_2 = 0.79061,\\
    & \varphi_1 = 2.09416, \quad \quad  \varphi_2 = 2.89674.
\end{align*}
A rigorous verification in Maple proves that this quintuple indeed corresponds to a solution of Rupert's problem for the Truncated icosidodecahedron. The visualization of this solution is presented in Figure~\ref{fig:truncated_icosidodecahedron} where the two projections of the polyhedron are plotted such that the black one lies inside the red one.
\end{proof}

We explained in the introduction that of the 13 Archimedean solids 8 were proven to be Rupert in~\cite{ChYaZa18} and an additional one in~\cite{Hoffmann19,Lavau19}. The theorem above raises this number up to 10, leaving only three Archimedean solids open: Snub cube, Rhombicosidodecahedron and Snub dodecahedron.

The same method and proof as above can be applied to the family of dual solids to the Archimedean ones, called the Catalan solids\footnote{For the coordinates of the Catalan solids we used we refer to the wonderful website \href{http://dmccooey.com/polyhedra/Catalan.html}{www.dmccooey.com/polyhedra/Catalan.html}.}. We obtain:
\begin{thm}\label{thm:catalan_rupert_solids}
The Rhombic dodecahedron, Triakis octahedron, Tetrakis hexahedron, Deltoidal icositetrahedron, Disdyakis dodecahedron, Rhombic triacontahedron, Triakis icosahedron, Pentakis dodecahedron and Disdyakis triacontahedron all have Rupert's property.
\end{thm}
\begin{proof}
The parameters for the solution of each solid are displayed in Table~\ref{tab:all_solutions}.
\end{proof}

Interestingly, this theorem shows that, similarly to Archimedean solids, 9 of the 13 Catalan solids admit Rupert's property. Except for the Triakis tetrahedron (Figure~\ref{fig:TriTetRID}, left), the remaining unresolved ones are precisely the dual polyhedra of the unsolved Archimedean solids. This raises the question on connectivity of the notions of duality and Rupert's property; we will state it precisely in \S\ref{sec:discussion}. 

As mentioned in the introduction, since the submission of this work, Fredriksson~\cite{Fredriksson22} could prove that the Triakis tetrahedron and Pentagonal icositetrahedron are Rupert.  
\\

In order to test the power of the presented algorithms, we ran our implementation on the family of 92 Johnson solids\footnote{Exact coordinates taken from  \href{http://dmccooey.com/polyhedra/Johnson.html}{www.dmccooey.com/polyhedra/Johnson.html}.}. We let the algorithm search for a solution for each polyhedron for at most an hour. The result is as follows.
\begin{thm}\label{thm:johnson_rupert}
Out of the 92 Johnson solids (at least) 82 admit Rupert's property. The remaining ones are: $J25$, $J45$, $J47$, $J71$, $J72$, $J73$, $J74$, $J75$, $J76$,  $J77$.
\end{thm}

Note that $J71$, $J72$, $J73$, $J74$, $J75$, $J76$,  $J77$ are all closely connected to the Rhombicosidodecahedron which we conjecture to be not Rupert (Conjecture~\ref{conj:RCD_non_Rupert}).

\subsection{Lower bounds on Nieuwland numbers} \label{sec:lowerNieuwland}

Running the algorithm from \S\ref{sec:findimproveNieuwland} for a few hours on the solved Platonic and Archimedean solids, we could significantly improve most of the previously known lower bounds for their Nieuwland numbers. Table~\ref{tab:Nieuwland_constants} summarizes these results. Like before, these numbers are found numerically in R and then verified rigorously in Maple. 

\begin{figure}[htb]
    \centering
    \begin{minipage}{0.49\textwidth}
        \includegraphics[width = 4.5cm]{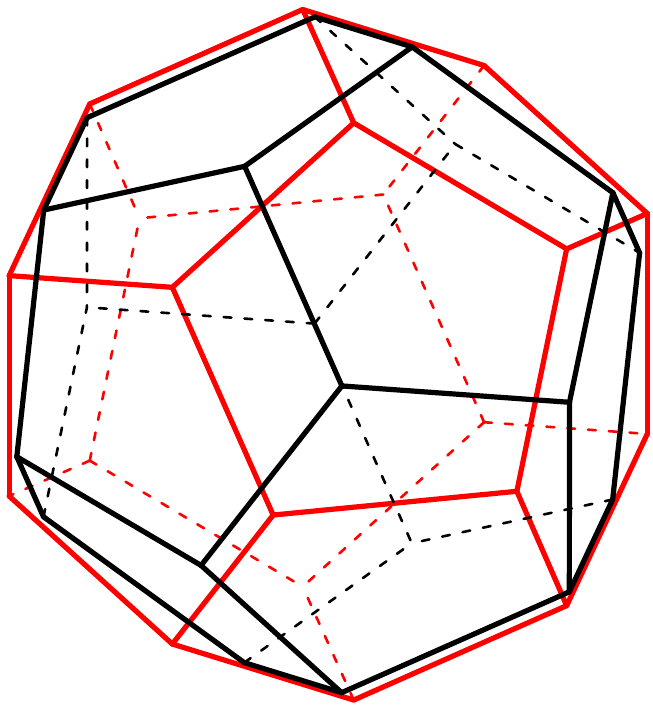}
    \end{minipage}
    \hfill
    \begin{minipage}{0.49\textwidth}
        \includegraphics[width = 5.3cm]{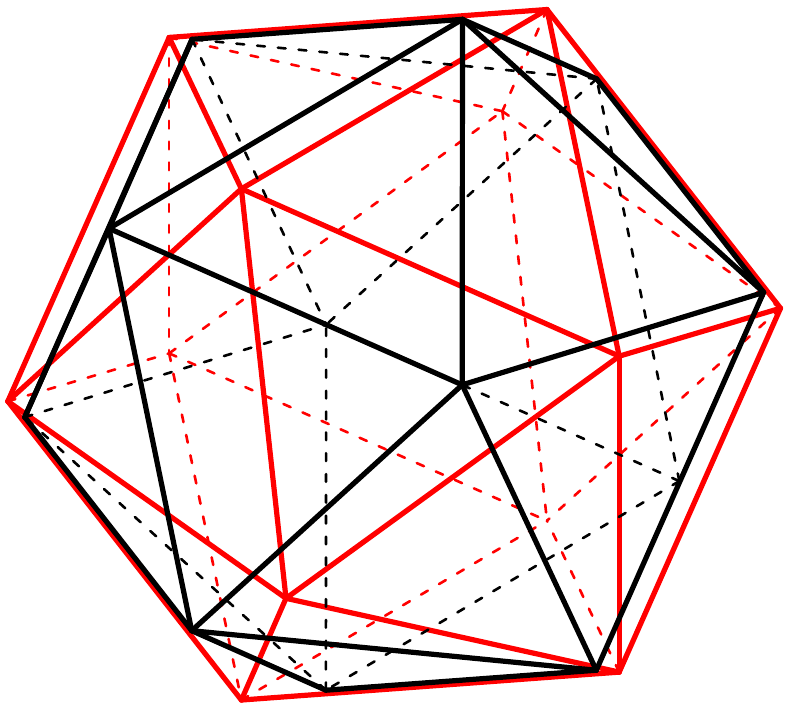}
    \end{minipage}
    \caption{Solution of the Rupert's problem for the Dodecahedron~(left) with Nieuwland number 1.010818 and for the Icosahedron~(right) with Nieuwland number 1.010805.}
    \label{fig:dodecahedronandicosahedron}
\end{figure}

For the Platonic solids Dodecahedron and Icosahedron we have found solutions with Nieuwland numbers 1.010818 and 1.010805 respectively. These figures are lower bounds for the Nieuwland numbers of these polyhedra. The numerical similarity of these numbers suggests that possibly they agree completely, like it is (conjecturally~\cite{JeWeYu17}) the case for the Cube and Octahedron. We address this question again in \S\ref{sec:discussion}. Figure~\ref{fig:dodecahedronandicosahedron} is a visualization of our solutions to Rupert's problem for the Dodecahedron and Icosahedron. In both cases we plot different projections of the solids in red and black such that the black projection lies inside the red one.

\begin{table}[t]
\begin{tabular}{l|c|c|r}
\hline
Name of solid               & Old best $\mu$& New best $\mu$ & Improvement \\ \hline \hline
Tetrahedron                 & $1.004\,235$  & $1.014\,473 $  &  3.42          \\ \hline
Cube                        & $1.060\,660$  & $1.060\,659$&  -- \\ \hline
Octahedron                  & $1.060\,660$  &$1.060\,640$& --  \\ \hline
Dodecahedron                & $1.005\,882$  &$1.010\,818$ &  1.84  \\ \hline
Icosahedron                 & $1.009\,107$  &$1.010\,805$ & 1.19 \\ \hline \hline
Truncated tetrahedron       & $>1$          & $ 1.014\,210$ &  --  \\ \hline
Cuboctahedron               & $1.014\,61$   &   $1.014\,571$ &  --  \\ \hline
Truncated cube              & $1.020\,36$   &  $1.030\,659$  &  1.51  \\ \hline
Truncated octahedron        & $1.008\,15$   &  $1.014\,602$ &   1.79  \\ \hline
Rhombicuboctahedron         & $1.006\,09$   & $1.012\,819$ & 2.10 \\ \hline
Truncated cuboctahedron     & $1.003\,70$   &  $1.006\,563$  &  1.77   \\ \hline
Snub cube                   & ---           &  ---      &   \\ \hline
Icosidodecahedron           & $1.000\,15$   &  $1.000\,878$ &  5.85  \\ \hline
Truncated dodecahedron      & $1.000\,14$   &  $1.001\,612$ & 11.51 \\ \hline
Truncated icosahedron       & $1.000\,04$   &  $1.001\,955$ &   48.88  \\ \hline
Rhombicosidodecahedron      & ---           & ---       &  -- \\ \hline
Truncated icosidodecahedron &  ---          & $1.002\,048$ &  -- \\ \hline
Snub dodecahedron           &  ---          &  ---      &   -- \\ \hline
\end{tabular}
\caption{Improved Nieuwland numbers for Platonic and Archimedean solids. The old best lower bounds for the Nieuwland numbers are taken from \cite{JeWeYu17} and \cite{ChYaZa18}. The improvement is calculated using $(\mu_\text{new}-1)/(\mu_{\text{old}}-1)$.\vspace{-0.75cm}}
\label{tab:Nieuwland_constants}
\end{table}

\subsection{Estimating Rupertness}
Recall from Definition~\ref{def:rupertness} that the Rupertness of a point symmetric polyhedron is the probability that a pair of uniformly random projections of it can be extended to a solution of Rupert's problem. Like we explained in \S\ref{sec:Rupertness}, we can estimate this probability by randomly drawing projections $M_{\theta_1,\varphi_1}(\textbf{P}),M_{\theta_2,\varphi_2}(\textbf{P})$ and then searching for $\alpha \in (0,\pi)$ such that $(R_{\alpha} \circ M_{\theta_1,\varphi_1}) (\textbf{P}) \subset M_{\theta_2,\varphi_2} (\textbf{P})^\circ$ holds. For each of the 14 point symmetric Platonic and Archimedean we drew at least 10 million pairs of random projections and for each pair decided on the existence of such an $\alpha \in (0,\pi)$. The quantities of corresponding solutions are summarized in Table~\ref{tab:all_solutions}. For example, the first row means that out of our $10^7$ random projections of the Cube precisely 657337 can be extended to a solution of Rupert's problem. This means that the Rupertness of the Cube is approximately 6.57\% and the 99.9\% confidence interval calculated with the Clopper-Pearson formula~(\ref{eq:Clopper-Pearson}) is $(0.0655 , 0.0659)$.

\begin{table}[h]
\centering
\begin{tabular}{l|c|c|c|c}
\hline
Name of solid             & $n$ & $k$ & $k/n $ & Confidence interval\\  & & & (in \%) &($\alpha = 99.9\%$)\\ \hline \hline
Cube                        & $10^7$ &657337& 6.57 &  $(0.0655 , 0.0659)$ \\ \hline
Octahedron                  & $10^7$  &1195417 & 11.95 & $(0.119, 0.120)$\\ \hline
Dodecahedron                & $10^7$  &230918& 2.31 & $(0.0230 , 0.0232)$ \\ \hline
Icosahedron                 & $10^7$  &295645& 2.96 & $(0.0294, 0.0297)$ \\ \hline \hline
Cuboctahedron               & $10^7$   & 390404& 3.90 &  $(0.0389, 0.0392)$ \\ \hline
Truncated cube              & $10^7$   &335602 & 3.36 &  $(0.0334, 0.0337)$ \\ \hline
Truncated octahedron        & $10^7$   &149188  & 1.49 &  $(0.0148, 0.0150)$\\ \hline
Rhombicuboctahedron         & $10^7$   &131176 & 1.31 &  $(0.0130, 0.0132)$\\ \hline
Truncated cuboctahedron     & $10^7$   &46044   & 0.460 & $(0.00455, 0.00466)$\\ \hline
Icosidodecahedron           & $10^7$   &40046 & 0.400 & $(0.00395, 0.00406)$ \\ \hline
Truncated dodecahedron      & $10^7$   &7583 & $0.0758$  & $(0.000736, 0.000781)$ \\ \hline
Truncated icosahedron       & $10^7$   &10813  & 0.108 & $(0.00105, 0.00111)$ \\ \hline
\textit{Rhombicosidodecahedron}      & $10^8$    & 0 & 0 & $[0,0.000000053)$ \\ \hline
Truncated &&&&\\
icosidodecahedron &  $10^7$     & 16394     & 0.164 & $(0.00161, 0.00167)$\\ \hline
\end{tabular}
\caption{Estimation of the Rupertness of point symmetric Platonic and Archimedean solids. The column $k$ says how many of the $n$ randomly chosen projections can be extended to solutions. $k/n$ is the estimate of the Rupertness and the last column is the $99.9\%$ confidence interval for it.}
\label{tab:Rupertness}
\end{table}

\begin{figure}[h!]
    \centering
        \includegraphics[width = 5cm]{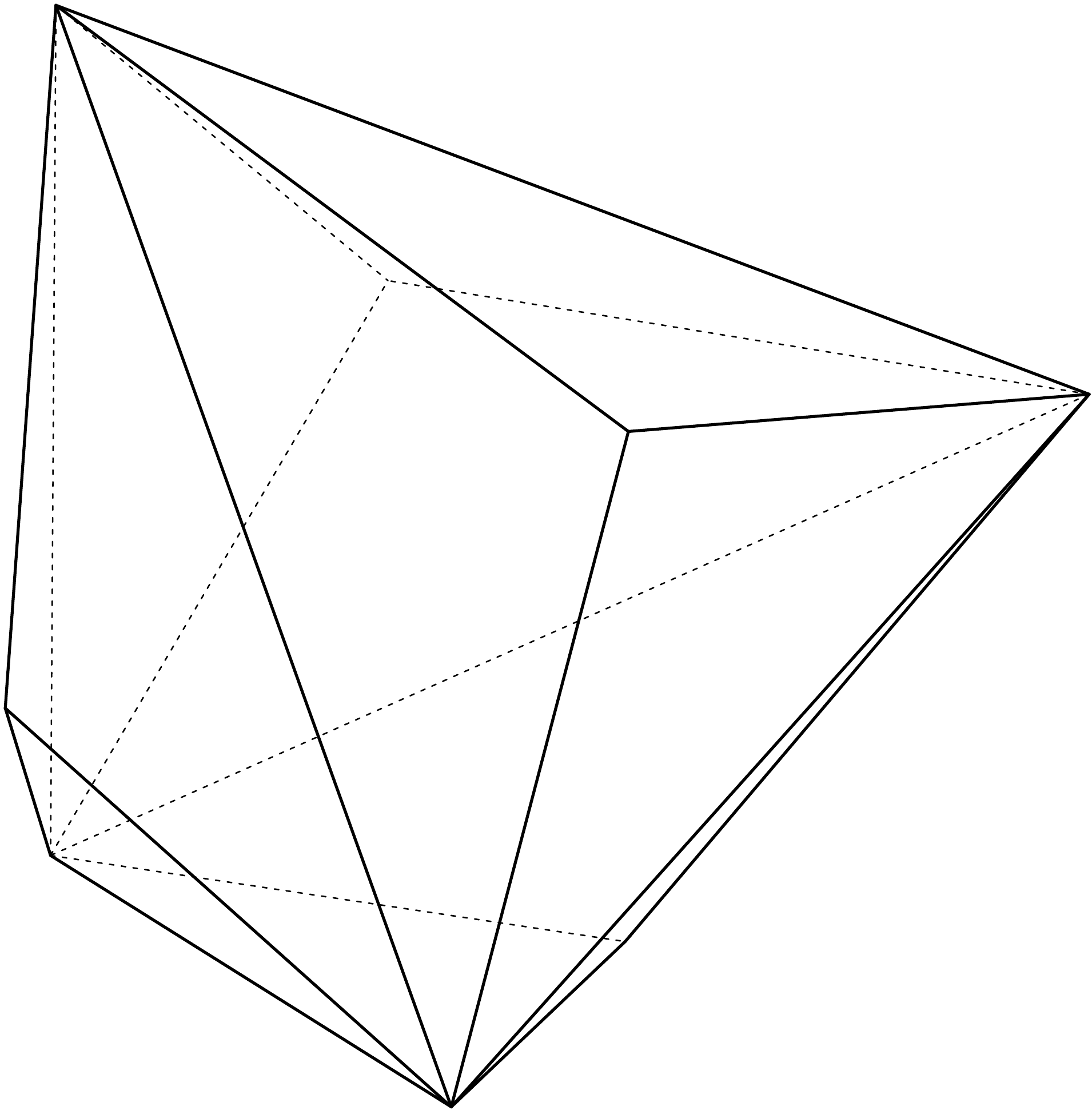}
        \quad
        \includegraphics[width=5cm]{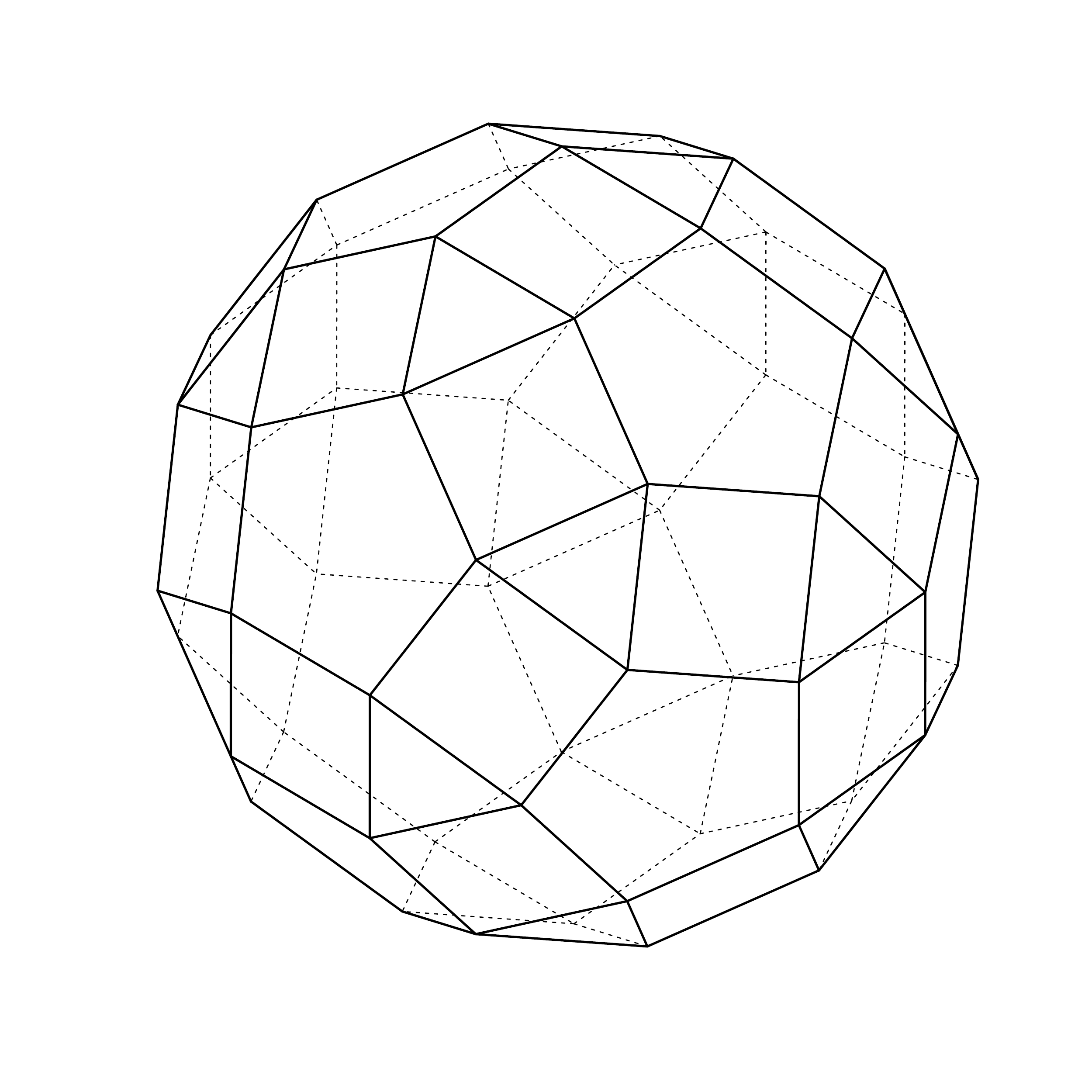}
    \caption{Triakis tetrahedron (left), Rhombicosidodecahedron (right)}
    \label{fig:TriTetRID}
\end{figure}

One notices immediately that the Rhombicosidodecahedron (Figure~\ref{fig:TriTetRID} right) is not only still unsolved regarding Rupert's property, since out of 100 million tries 0 could have been extended to a solution, but also that its Rupertness is (with confidence of 99.9\%) significantly lower than the Rupertness of any other point symmetric Platonic or Archimedean solid. In fact, with probability 99.9\%, the Rupertness of the RID is less than 1/10000 of the Rupertness of the Truncated dodecahedron, the one with the smallest figure. Based on Table~\ref{tab:Rupertness} we state the following surprising conjecture which contradicts Conjecture~\ref{conj:every_poly_rupert} taken from \cite[Open~problem,~Conjecture,~p. 503]{ChYaZa18}. 
\begin{conj}\label{conj:RCD_non_Rupert}
The Rhombicosidodecahedron does not have Rupert's property.
\end{conj}
\bigskip

A natural attempt to prove Conjecture~\ref{conj:RCD_non_Rupert} would be to employ the deterministic algorithm in Theorem~\ref{thm:deterministic}, or rather its extension to polyhedra with coordinates given by algebraic numbers (see \S\ref{sec:deterministic}). Like we already explained in the remark at the end of \S\ref{sec:deterministic}, the bound $n!$ is very pessimistic also in this case. To be precise, we are confident that it should not be difficult to prove that (accounting for symmetries) there are not more than {50} possible silhouettes to consider for the Rhombicosidodecahedron. Since the RID has 60 vertices, it follows that we would need to prove emptiness of {50} semi-algebraic sets defined by at most $60^2 = 3600$ polynomial inequalities in $7-2+1=6$ variables (we can set $x=y=0$ but we need a variable for the golden ratio) of total degree of at most 22. Unfortunately, it seems that these numbers are too big for current algorithms and implementations: in order to have a chance for termination in reasonable time, we would need to reduce the number of inequalities to below 20. Therefore, Conjecture~\ref{conj:RCD_non_Rupert} is still open.

Initially we were quite skeptical that the other unsolved Archimedean solids (Snub cube and Snub Dodecahedron) as well as for the four unsolved Catalan solids (numbers 19, 25, 29, 31 in Table~\ref{tab:all_solutions}) and the 10 open Johnson solids (see Theorem~\ref{thm:johnson_rupert}) admit Rupert's property. For these solids we did not estimate the Rupertness and hence have no statistical evidence; so we concluded that it is very much possible that one should just execute the algorithms for a longer time in order to find a solution. Indeed, after the submission of our work, Fredriksson~\cite{Fredriksson22} was able to improve on our methods and show Rupert's property for the Catalan solids 19 and 25 in Table~\ref{tab:all_solutions}, as well as for the Johnson solids J25, J45, J47, J71 and J76. We concentrated our search on the RID, since it is the smallest point symmetric solid for which we could not find a solution to Rupert's problem.

\noindent

\pagebreak

\subsection{Concluding remarks and future work} \label{sec:discussion}
One may notice a surprising fact in Theorem~\ref{thm:catalan_rupert_solids}: A point symmetric Archimedean solid is proven to be Rupert if and only if its dual solid is. While this is only a small indication for the connectivity of duality and Rupert's property, Table~\ref{tab:Nieuwland_constants} provides more evidence: the Cube and the Octahedron are conjectured to have the same Nieuwland number and the same seems to hold for the other pair of dual Platonic solids: the Dodecahedron and Icosahedron. Based on these observations we formulate natural and interesting but apparently not easy-to-answer questions:
\begin{enumerate}
    \item Is a point symmetric Archimedean solid Rupert if and only if its dual Catalan solid has Rupert's property?
    \item Do dual Platonic solids have the same Nieuwland number? If so, is there a geometric reason for this?
    \item What are the exact Nieuwland numbers of the Dodecahedron and Icosahedron? Do they also admit simple algebraic expressions like the (conjectured) $3\sqrt{2}/4$ for the Cube and Octahedron?\footnote{In~\cite{StYu22} we conjecture that the minimal polynomial for these numbers is given by $2025x^8 - 11970x^6 + 17009x^4 - 9000x^2 + 2000$.} \label{question:3}
\end{enumerate}
If a solution to Rupert's problem of a Platonic or Archimedean solid is given in $\R^3$ by $\textbf{P}$ and its copy $\textbf{Q}$, one can look at the duals of both polyhedra. It is quite intriguing that it seems that the dual of an ``optimal'' solution (i.e. one with highest Nieuwland number) of a Platonic solid yields an ``optimal'' solution for the dual solid. However, we could not find a (geometric) explanation for this. Moreover, the dual of \textit{some} solution of a Platonic or Archimedean solid is not necessarily a solution at all. 

As already mentioned, Conjecture~\ref{conj:RCD_non_Rupert} contradicts current beliefs on Rupert's property for polyhedra, but at the same time we have statistical reasons to believe in our conjecture. Assuming its validity, further natural questions are:
\begin{enumerate}
\setcounter{enumi}{3}
    \item What distinguishes the Rhombicosidodecahedron from other (point symmetric) Archimedean solids and prevents this polyhedron to have Rupert's property? Is there an easy criterion for Rupert polyhedra? 
    \item How can one prove Conjecture~\ref{conj:RCD_non_Rupert}? Are the remaining Archimedean and Catalan solids (12, 16, 18, 19, 25, 29, 31 in Table~\ref{tab:all_solutions}) Rupert?\footnote{The Catalan solids 19 and 25 are now resolved in~\cite{Fredriksson22}.}
\end{enumerate}

\noindent
\textbf{Acknowledgments}\\
We would like to thank Alin Bostan for his careful and supportive reading of the first version of this manuscript, Herwig Hauser for his constantly encouraging advices, and Mohab Safey El Din for his great help and expertise on algorithms for semi-algebraic sets. The authors are also grateful to David I. McCooey for creating the wonderful website \href{http://dmccooey.com/polyhedra/}{www.dmccooey.com/polyhedra} which contains coordinates, visualizations and information for hundreds of interesting polyhedra. Finally, the authors thank the anonymous referees for their helpful comments.

The second author was financially supported by the DOC fellowship of the \href{https://www.oeaw.ac.at/en/}{ÖAW} (26101), the WTZ collaboration project of the \href{https://oead.at/en/}{OeAD} (FR 09/2021) and the {\href{https://specfun.inria.fr/chyzak/DeRerumNatura/}{DeRerumNatura}} project ANR-19-CE40-0018.

\section{Appendix} \label{sec:appendix}

In the appendix we most importantly present Table~\ref{tab:all_solutions} which summarizes our solutions to Rupert's problem for all Platonic, 10 Archimedean and 9 Catalan solids. According to Proposition~\ref{prop:Rupert_7parameters}, any solution can be encoded by seven parameters $x,y,\alpha,\theta_1,\theta_2,\varphi_1,\varphi_2$. So for each solved polyhedron we provide these numbers in the corresponding columns. Proposition~\ref{prop:pointsym} implies that if a polyhedron is point symmetric, one can choose $x=y=0$, so in these cases $x$ and $y$ are zero. The right column of Table~\ref{tab:all_solutions} shows the Nieuwland number of the solution.

Finally, Table \ref{tab:Coordinates} incorporates the exact coordinates we used for the Platonic and Archimedean solids. The coordinates for Catalan and Johnson solids can be found at \href{https://github.com/Vog0/RupertProblem/}{www.github.com/Vog0/RupertProblem} and are taken, as mentioned, from the website \href{http://dmccooey.com/polyhedra}{www.dmccooey.com/polyhedra}. The first link also contains the source code in R and Maple we used to find and then verify solutions.

{\small
\begin{table}[H]
\vspace*{0cm}
\begin{adjustwidth}{-1.3cm}{}
\centering
\Rotatebox{90}{
    \begin{tabular}{r|l|l|l|l|l|l|l|l|r}
    Nr. & Name of solid & $x$ & $y$ & $\alpha$ & $\theta_1$ & $\varphi_1$ & $\theta_2$ & $\varphi_2$ & $\mu(v,\textbf{P})$ \\
    \hline \hline
        1. & Tetrahedron & 0.1788244 & -0.0976062 & 1.0372426 & 5.3278439 & 1.5713832 & 3.9444529 & 0.9501339 & 1.014473 \\ \hline
        2. & Cube & 0 & 0 & 2.4840821 & 1.9060829 & 3.1415929 & 5.8188256 & 2.3004443 & 1.060659 \\ \hline
        3. & Octahedron & 0 & 0 & 3.1415873 & 5.4977985 & 1.9105975 & 6.2808288 & 1.5701448 & 1.060640 \\ \hline
        4. & Dodecahedron & 0 & 0 & 1.0378047 & 0.8553414 & 2.108091 & 4.918788 & 2.0545287 & 1.010818 \\ \hline
        5. & Icosahedron & 0 & 0 & 2.7276836 & 2.7732324 & 2.6181502 & 2.3091726 & 2.2712915 & 1.010805 \\ \hline \hline
        6. & Truncated tetrahedron & 0.160858 & -0.164724 & 4.7775741 & 6.2831072 & 0.7854425 & 2.0992734 & 1.3849498 & 1.014210 \\ \hline
        7. & Cuboctahedron & 0 & 0 & 3.1386793 & 2.5259348 & 1.5710827 & 0.7902177 & 0.9351593 & 1.014571 \\ \hline
        8. & Truncated cube & 0 & 0 & 2.298646 & 4.3427928 & 3.1415862 & 2.089632 & 2.2876946 & 1.030659 \\ \hline
        9. & Truncated octahedron & 0 & 0 & 1.5690349 & 3.1415601 & 0.785367 & 5.3243536 & 2.0933886 & 1.014602 \\ \hline
        10. & Rhombicuboctahedron & 0 & 0 & 0.017061 & 2.9503929 & 3.1415921 & 4.1693802 & 0.636201 & 1.012819 \\ \hline
        11. & Truncated cuboctahedron & 0 & 0 & 0.2396229 & 3.1416249 & 0.785486 & 4.4525352 & 0.429099 & 1.006563 \\ \hline
        12. & Snub cube & -- & -- & -- & -- & -- & -- & -- & -- \\ \hline
        13. & Icosidodecahedron & 0 & 0 & 1.578603 & 2.7736451 & 0.7120286 & 4.7086522 & 2.1263666 & 1.000878 \\ \hline
        14. & Truncated dodecahedron & 0 & 0 & 2.2092757 & 4.3599229 & 1.5508055 & 1.6477247 & 1.0979977 & 1.001612 \\ \hline
        15. & Truncated icosahedron & 0 & 0 & 0.9547212 & 4.7124428 & 1.470154 & 0.8649729 & 2.0954566 & 1.001955 \\ \hline
        16. & Rhombicosidodecahedron & -- & -- & -- & -- & -- & -- & -- & -- \\ \hline
        17. & Truncated icosidodecahedron & 0 & 0 & 0.4358364 & 2.7768504 & 2.0941596 & 0.79061 & 2.8967442 & 1.002048 \\ \hline
        18. & Snub dodecahedron & -- & -- & -- & -- & -- & -- & -- & -- \\ \hline \hline
        19. & Triakis tetrahedron & -- & -- & -- & -- & -- & -- & -- & -- \\ \hline
        20. & Rhombic dodecahedron & 0 & 0 & 0.2389694 & 3.926939 & 0.9553557 & 5.171164 & 1.3442843 & 1.027201 \\ \hline
        21. & Triakis octahedron & 0 & 0 & 0.3562255 & 5.7674031 & 2.2867379 & 0.0005374 & 1.5665899 & 1.030648 \\ \hline
        22. & Tetrakis hexahedron & 0 & 0 & 0.1945682 & 3.4241341 & 1.1711373 & 0.0040963 & 2.3603178 & 1.009632 \\ \hline
        23. & Deltoidal icositetrahedron & 0 & 0 & 0.6277374 & 0.6012867 & 1.4476059 & 6.1255227 & 3.1382821 & 1.007632 \\ \hline
        24. & Disdyakis dodecahedron & 0 & 0 & 0.1178211 & 6.1466092 & 2.5957828 & 1.5695218 & 0.7842378 & 1.002500 \\ \hline
        25. & Pentagonal icositetrahedron & -- & -- & -- & -- & -- & -- & -- & -- \\ \hline
        26. & Rhombic triacontahedron & 0 & 0 & 0.231712 & 2.84e-05 & 0.5535717 & 1.9227518 & 2.1379305 & 1.007037 \\ \hline
        27. & Triakis icosahedron & 0 & 0 & 2.5481489 & 3.3133906 & 0.4995076 & 2.3963212 & 2.1824603 & 1.001304 \\ \hline
        28. & Pentakis dodecahedron & 0 & 0 & 3.1547479 & 5.4202246 & 2.1024926 & 4.2553188 & 2.4568193 & 1.001845 \\ \hline
        29. & Deltoidal hexecontahedron & -- & -- & -- & -- & -- & -- & -- & -- \\ \hline
        30. & Disdyakis triacontahedron & 0 & 0 & 2.5886126 & 4.2871288 & 0.7860227 & 5.917639 & 2.107937 & 1.000210 \\ \hline
        31. & Pentagonal hexecontahedron & -- & -- & -- & -- & -- & -- & -- & -- \\ \hline
    \end{tabular}
}
\caption{Solutions to Rupert's problem for Platonic, Archimedean and Catalan solids.}
\label{tab:all_solutions}
\end{adjustwidth}
\end{table}
}

{\small
\begin{table}[H]
\begin{adjustwidth}{-1cm}{}
\centering
\begin{tabular}{l|l}
\hline
Name of solid  & Coordinates \\ \hline \hline
1. Tetrahedron                  & $(\pm1,\pm1,\pm1)$ with an even number of ``$-$" signs               \\ \hline
2. Cube                         & $(\pm1,\pm1,\pm1)$                 \\ \hline
3. Octahedron                   & all permutations of $(0,0,\pm1)$                 \\ \hline
4. Dodecahedron                 & $(\pm1,\pm1,\pm1)$ and all even permutations of $(0,\pm \Phi^{-1}, \pm \Phi)$                \\ \hline
5. Icosahedron                  & even permutations of $(0,\pm\Phi,\pm1)$                \\ \hline \hline
6. Truncated tetrahedron        &  all permutations of $(\pm1,\pm1,\pm3)$ with an even number \\
                                &of ``$-$" signs \\ \hline
7. Cuboctahedron                & all permutations of $(\pm1,\pm1,0)$\\ \hline
8. Truncated cube               & all permutations of $(\pm1,\pm1,\pm (\sqrt{2}-1))$ \\ \hline
9. Truncated octahedron         &  all permutations of $(0,\pm1,\pm 2)$\\ \hline
10. Rhombicuboctahedron         & all permutations of $(\pm1,\pm1,\pm (1+\sqrt{2}))$\\ \hline
11. Truncated cuboctahedron     & all permutations of $(\pm1,\pm(1+\sqrt{2}),\pm (1+2\sqrt{2}))$\\ \hline
12. Snub cube                   & all even permutations of $(\pm1, \pm1/t,\pm t)$
with an even number\\
&of plus signs and all odd permutations with an odd number\\
&of plus signs. $t$ is the tribonacci constant             \\ \hline
13. Icosidodecahedron           & all permutations of $(0,0,\pm \Phi)$ and all even permutations \\
                                &of $(\pm\frac{1}{2},\pm\frac{\Phi}{2},\pm\frac{\Phi^2}{2})$\\ \hline
14. Truncated dodecahedron      & all even permutations of $(0,\pm1/\Phi,\pm(2+\Phi)), (\pm\frac{1}{\Phi},\pm\Phi,\pm2\Phi)$\\
                                &and $(\pm\Phi,\pm2,\pm(\Phi+1))$   \\ \hline
15. Truncated icosahedron       & all odd permutations of $(0,\pm1,\pm3\Phi), (\pm1,\pm(2+\Phi),\pm2\Phi)$ \\
                                & and $(\pm\Phi,\pm2,\pm(2\Phi+1))$     \\ \hline
16. Rhombicosidodecahedron      & all even permutations of $(\pm1,\pm1,\pm\Phi^3), (\pm\Phi^2,\pm\Phi,\pm2\Phi)$\\
                                &and $(\pm(2+\Phi),0,\pm\Phi^2)$     \\ \hline
17. Truncated icosidodecahedron & all even permutations of $(\pm\frac{1}{\Phi},\pm\frac{1}{\Phi},\pm(3+\Phi)),$\\
&$(\pm\frac{2}{\Phi},\pm\Phi,\pm(1+2\Phi)),(\pm\frac{1}{\Phi},\pm\Phi^2,\pm(-1+3\Phi)),$\\
&$(\pm(2\Phi-1),\pm2,\pm  (2+\Phi))$ and $(\pm\Phi,\pm3,\pm2\Phi)$     \\ \hline
18. Snub dodecahedron           & all even permutations of $(2\alpha,2,2\beta)$,\\
&$(\alpha+\frac{\beta}{\Phi}+\Phi,-\alpha\Phi+\beta+\frac{1}{\Phi},\frac{\alpha}{\Phi}+\beta\Phi-1)$,\\
&$(\alpha+\frac{\beta}{\Phi}-\Phi,\alpha\Phi-\beta+\frac{1}{\Phi},\frac{\alpha}{\Phi}+\beta\Phi+1)$,\\
&$(-\frac{\alpha}{\Phi}+\beta\Phi+1,-\alpha+\frac{\beta}{\Phi}-\Phi,\alpha\Phi+\beta-\frac{1}{\Phi})$,\\
&$(-\frac{\alpha}{\Phi}+\beta\Phi-1,\alpha-\frac{\beta}{\Phi}-\Phi,\alpha\Phi+\beta+\frac{1}{\Phi})$\\
&with an odd number of sign changes of the coordinates, where\\
&$\xi = \sqrt[3]{\frac{\Phi}{2} + \frac{1}{2}\sqrt{\Phi - \frac{5}{27}}} + \sqrt[3]{\frac{\Phi}{2} - \frac{1}{2}\sqrt{\Phi - \frac{5}{27}}}$, $\alpha=\xi-1/\xi$ and\\
&$\beta=\xi\Phi+\Phi^2+\Phi/\xi$
                  \\ \hline
\end{tabular}
\caption{Used coordinates of all Platonic and Archimedean solids, as used in the Maple Package \href{https://maplesoft.com/support/help/maple/view.aspx?path=geom3d}{geom3d} (for verification), except Snub Cube and Snub Dodecahedron, which are not needed for our results. $\Phi = (\sqrt{5}+1)/2 \approx 1.62$ is the golden ratio.}
\label{tab:Coordinates}
\end{adjustwidth}
\end{table}
}

\bibliographystyle{alphaabbr}
\bibliography{bib}

\vspace{1cm}

\noindent

\end{document}